\documentclass[oneside, 12pt]{amsart}
\usepackage{amscd, amssymb, amsmath, mathrsfs}
\usepackage[english]{babel}
\usepackage{booktabs}
\usepackage{tikz-cd}
\usepackage{url}
\usepackage[pdftex, colorlinks=true,  citecolor=blue, linkcolor=blue, linktocpage=true]{hyperref} 

\newcommand\mylabel[1]{\label{#1}\marginpar{\vspace{-1ex}\medskip\medskip\footnotesize \tt #1}}
\renewcommand\mylabel[1]{\label{#1}}
\newcommand{\mydate}{
\number\day\space
\ifcase\month \or January\or February\or March\or April\or May\or June\or July\or August\or September\or October\or November\or December\fi 
\space\number\year}

\DeclareUrlCommand\arXiv{\urlstyle{same}}

\newtheorem{theorem}{Theorem}[section]
\newtheorem{maintheorem}{Theorem}
\newtheorem{lemma}[theorem]{Lemma}
\newtheorem{proposition}[theorem]{Proposition}
\newtheorem{corollary}[theorem]{Corollary}

\theoremstyle{definition}
\newtheorem{definition}[theorem]{Definition}

\newtheorem*{acknowledgement}{Acknowledgement}

\theoremstyle{remark}

\newcommand{\ZZ}{\mathbb{Z}}
\newcommand{\QQ}{\mathbb{Q}}

\newcommand{\PP}{\mathbb{P}}

\newcommand{\ideala}{\mathfrak{a}}
\newcommand{\idealb}{\mathfrak{b}}
\newcommand{\idealc}{\mathfrak{c}}

\newcommand{\shE}{\mathscr{E}}
\newcommand{\shF}{\mathscr{F}}

\newcommand{\shH}{\mathscr{H}}

\newcommand{\shI}{\mathscr{I}}

\newcommand{\shN}{\mathscr{N}}
\newcommand{\shL}{\mathscr{L}}

\newcommand{\alg}{\text{\rm alg}}

\newcommand{\Bl}{\operatorname{Bl}}

\newcommand{\edim}{\operatorname{edim}}

\newcommand{\Exc}{\operatorname{Exc}}
\newcommand{\Ext}{\operatorname{Ext}}

\newcommand{\GL}{\operatorname{GL}}

\newcommand{\Hom}{\operatorname{Hom}}

\newcommand{\id}{{\operatorname{id}}}

\newcommand{\length}{\operatorname{length}}

\newcommand{\lra}{\longrightarrow}

\newcommand{\maxid}{\mathfrak{m}}

\newcommand{\Num}{\operatorname{Num}}

\renewcommand{\O}{\mathscr{O}}

\newcommand{\pdeg}{\operatorname{pdeg}}

\newcommand{\Pic}{\operatorname{Pic}}
\newcommand{\PGL}{\operatorname{PGL}}

\newcommand{\pr}{\operatorname{pr}}
\newcommand{\Proj}{\operatorname{Proj}}

\newcommand{\quadand}{\quad\text{and}\quad}

\newcommand{\ra}{\rightarrow}

\newcommand{\red}{{\operatorname{red}}}
\newcommand{\Reg}{\operatorname{Reg}}

\newcommand{\Sing}{\operatorname{Sing}}

\newcommand{\Spec}{\operatorname{Spec}}

\newcommand{\uHom}{\underline{\operatorname{Hom}}}

\begin{document}

\title[The Iskovskih  Theorem]
      {The Iskovskih  Theorem for regular surfaces over imperfect fields}

\author[Andrea Fanelli]{Andrea Fanelli}
\address{Univ. Bordeaux, CNRS, Bordeaux INP, IMB, UMR 5251, F-33400 Talence, France.}
\curraddr{}
\email{andrea.fanelli@math.u-bordeaux.fr}

\author[Stefan Schr\"oer]{Stefan Schr\"oer}
\address{Heinrich Heine University D\"usseldorf, Faculty of Mathematics and Natural Sciences, Mathematical Institute, 40204 D\"usseldorf, Germany}
\curraddr{}
\email{schroeer@math.uni-duesseldorf.de}


\dedicatory{13 March 2025}

\begin{abstract}
We generalize Iskovskih's theorem about   surfaces   without irregularity and bigenus from the smooth case  
to regular surfaces over  arbitrary  fields,  with special focus on the case of imperfect fields. 
This includes surfaces that are geometrically non-normal or geometrically non-reduced. 
Here the usual approach of Galois descent breaks down, and one
relies entirely on the scheme theory over the ground field. Moreover, the  degrees of   closed points can  be larger than
expected, and certain curves might have purely inseparable constant field extension. To deal with the latter
we establish a general theory for  inseparable pencils, which is of independent interest.
A crucial case  not present in the classical proof for   Iskovskih's theorem leads to non-normal quartic surfaces that are singular
along a twisted cubic, or more exotic space curves of degree three. 
\end{abstract}

\maketitle
\tableofcontents

\section*{Introduction}
\mylabel{Introduction}

A cornerstone in the Enriques Classification of algebraic surfaces is  
\emph{Castelnuovo's Rationality Criterion}: A smooth surface $S$ over an algebraically closed ground field $k$
is birational to $\PP^2$ if and only if   $h^1(\O_S)=0$ and $h^0(\omega_S^{\otimes 2})=0$.
Note that $H^1(S,\O_S)$ is the Lie algebra for the Picard scheme $\Pic_{S/k}$, whose dimension
$h^1(\O_S)$ is classically called the \emph{irregularity}. The invariants $h^0(\omega^{\otimes n})$ stemming
from the dualizing sheaf $\omega_S=\det(\Omega^1_{S/k})$ are the \emph{plurigenera},
and $h^0(\omega_S^{\otimes 2})$ might be  called   \emph{bigenus}.

The above geometric statement can be seen as a consequence of \emph{Iskovskih's Theorem} \cite{Iskovskih 1980}, which 
is of arithmetic nature:
Let  $S$ be a \emph{minimal  smooth} surface  with $h^1(\O_S)=h^2(\omega_S^{\otimes 2})=0$, over an \emph{arbitrary} ground field $k$.
Then $S$ is isomorphic to the projective plane  or a quadric surface,  or there
is a fibration  where the base and the generic fiber are Brauer--Severi curves, 
or the dualizing sheaf $\omega_S$ generates the Picard group $\Pic(S)$.
This  builds on previous work of Manin \cite{Manin 1966}  over  perfect   fields $k$.

The goal of this paper is to generalize Iskovskih's Theorem, to allow \emph{minimal regular} surface  over  
arbitrary ground fields, with special focus on the case of  imperfect ground fields, similar in style to 
the investigations in \cite{Fanelli; Schroeer 2020a}.
Note that this  includes surfaces that are \emph{geometrically non-normal} or even \emph{geometrically non-reduced}. 
Our main result is:

\begin{maintheorem}
(see Thm.\ \ref{trichotomy})
Let $X$ be a minimal regular surface over an arbitrary field $F$, with numerical invariants
$h^0(\O_X)=1$ and  $h^1(\O_X)=h^0(\omega_X^{\otimes 2})=0$. Then one of the following holds:
\begin{enumerate}
\item The surface $X$ is isomorphic to the projective plane $\PP^2$,  or isomorphic to a   quadric surface  in $\PP^3$.
\item There is a morphism $f:X\ra B$ with $f_*(\O_X)=\O_B$ such that the base and the generic fibers are regular genus-zero curves.
\item The dualizing sheaf $\omega_X$ generates the Picard group $\Pic(X)$.
\end{enumerate}
\end{maintheorem}

In the first case of this trichotomy, the surface $X$ can be expressed in terms of a single   equation.
In the second case, the geometry is reduced to dimension one. In marked contrast, it is very difficult to say more
about the third case 
(compare \cite{Bernasconi; Fanelli; Schneider; Zimmermann 2024}, Section 4).
 
In recent years,    regular del Pezzo surface over imperfect have received increased attention, for their own sake, and  
because they appear as generic fibers
in the outputs of the   Minimal Model Program for algebraic varieties   in characteristic $p>0$,
for example in the work of  Bernasconi, Ji, Maddock, Martin, Patakfalvi, Tanaka, and ourselves
(\cite{Schroeer 2007a}, \cite{Maddock 2016} \cite{Tanaka 2018} \cite{Fanelli; Schroeer 2020a},  \cite{Bernasconi; Tanaka 2022},
\cite{Tanaka 2020}, \cite{Ji; Waldron 2021}, \cite{Patakfalvi; Waldron 2022}, \cite{Tanaka 2024},  \cite{Bernasconi; Tanaka 2024}, \cite{Bernasconi; Martin 2024}).
 
Koll\'ar, Smith and Corti  gave a highly readable presentation  of Iskovskih's arguments (\cite{Kollar; Smith; Corti 2004}, Chapter 3), and we 
follow their line of reasoning: Assuming that the dualizing shed does not generate the Picard group,
one produces  a   curve $C\subset X$ with $h^1(\O_C)=0$ and $C^2\geq 0$, such that $C$ and all linearly equivalent curves
are integral. \emph{Now two major new issues arise : The field of constants $E=H^0(C,\O_C)$ could be a purely inseparable extension,
and the self-intersection number $C^2\geq 0$ is potentially larger that in the classical situation.}
The first crucial step is to establish $h^0(\O_C)=1$. Having this, we infer  $C^2\mid 4$ or   $C^2=0$,
and the cases $0\leq C^2\leq 3$ than correspond to (i)--(iii) in Theorem A. The second crucial step is to establish that $C^2=4$ actually does not appear.
 
To tackle the problem of constant field extensions, we develop a \emph{general theory of purely inseparable pencils}.
Suppose $X$ is a proper normal scheme with $h^0(\O_X)=1$. Let   $D_0,D_1$ be two effective Cartier divisor stemming from global sections of
some invertible sheaf $\shL$, such that $Z=D_0\cap D_1$ is reduced and of codimension two.
The blowing-up  $V=\Bl_Z(X)$ comes with a fibration $h:V\ra\PP^1$, and we write $B=\Spec h_*(\O_V)$ for
the Stein factorization.  Our second main result reveals that   the prime $p=2$ plays a particular role:

\begin{maintheorem}
(see Thm.\ \ref{pencils characteristic two})
In the above setting, suppose also that the cohomology group $H^1(X,\O_X)$ vanishes,
and that for each rational point $t\in\PP^1$, the resulting effective Cartier divisor $D_t\subset X$ is reduced and geometrically connected.
The the following holds:
\begin{enumerate}
\item The Stein factorization $B$ is a regular genus-zero curve.
\item The map $g:B\ra\PP^1$ is a universal homeomorphism, with $\deg(g)\mid 2$.
\item For each non-zero $s\in H^0(X,\shL)$ the resulting $D\subset X$ has  $h^0(\O_D)=\deg(g)$.
\item If $\deg(g)=2$ then the ground field $F$ is imperfect of characteristic $p=2$, and $B$ is a twisted line or a twisted  ribbon.
\end{enumerate}
\end{maintheorem}

Here \emph{twisted lines} and \emph{twisted ribbons}  designate regular genus-zero  curves  without rational points that become, 
after   ground field extension,  isomorphic $\PP^1$ and the infinitesimal thickening
$\PP^1\oplus\O_{\PP^1}(-1)$, respectively. 

Now back to our minimal regular surface $X$ without irregularity and bigenus. To tackle the problem
of self-intersection numbers, we now assume that there is an integral curve $C\subset X$ with $h^1(\O_C)=0$ and $C^2=4$.
The following \emph{bend-and-break type statement} is our third main result:

\begin{maintheorem}
(see Thm.\ \ref{pencils characteristic two})
In the above situation, the curve $C\subset X$ is linearly equivalent to a curve $C'$ that is not integral.
\end{maintheorem}

For this we examine  four-dimensional linear systems for $\shL=\O_X(C)$ to produce  a morphism $f:X\ra \PP^3$
with $\shL=f^*\O_{\PP^3}(1)$. Using the theory of quadric hypersurfaces, it is not difficult to reduce to the case
that the image is a \emph{non-normal quartic surface}  $V\subset \PP^3$, and $f:X\ra V$ is the minimal resolution of singularities. 
Note that singular quartic surfaces form a classical evergreen topic,  as witnessed by the monographs of Hudson \cite{Hudson 1905}
and  Jessop \cite{Jessop 1916} from the beginning of the 19th century. For modern research, we mention
the work of Urabe (\cite{Urabe 1984}, \cite{Urabe 1986}, \cite{Urabe 1987}) and Catanese and Sch\"utt
(\cite{Catanese 2022}, \cite{Catanese; Schuett 2022}, \cite{Catanese; Schuett 2023}).

It turns out that the branch curve $B$ for the normalization map $g:Y\ra V$ is a space curve of degree three inside $\PP^3$.
The relevant case is when it is integral, and there are two possibilities:
Either $B$ is a \emph{twisted cubic}, a classical terminology that designates
a copy of the projective line embedded into the projective  three-space via $\O_{\PP^1}(3)$.
Or it is  a \emph{exotic cubic}, our ad hoc terminology to designates   non-Gorenstein genus-zero curves
arising as denormalizations of  the projective line over a cubic field extension.
The latter case can be treated via  the  conductor squares for the normalization of the  \emph{ramification curve} $R=g^{-1}(B)$,
involving commutative algebra for local Artin rings.

For the twisted cubics $B$, an entirely different strategy is required, and we  exploit the amazing  fact
that $\Bl_B(\PP^3)$ is in an unexpected way a $\PP^1$-bundle     over $\PP^2$, 
compare the work of Blanc and Lamy \cite{Blanc; Lamy 2012}, Ray \cite{Ray 2020}, and Sarkar \cite{Sarkar 2024}.
This ensures that  the induced map $\Bl_B(V)\ra\PP^2$ factors over a regular genus-zero  curve $V_+(\Phi)$ 
defined by some quadratic polynomial $\Phi\in F[X_0,X_1,X_2]$,
and $X$ becomes a ruled surface over such curves.

To unravel the geometry of the situation, we are now forced to classify the   
locally free sheaves on such twisted ribbons $B$. This can be seen as a variant
of \emph{Grothendieck's Splitting Theorem} for the projective line (see \cite{Grothendieck 1957} and  \cite{Hazewinkel; Martin 1982}),
and gives our  fourth main result:

\begin{maintheorem}
(See Thm.\ \ref{indecomposable on genus-zero})
Up to isomorphism, the indecomposable locally free sheaves on twisted ribbons $B$ are the 
$\omega_B^{\otimes a}$ and $\shF_B\otimes\omega_B^{\otimes b}$, with exponents $a,b\in\ZZ$.
\end{maintheorem}

Here  $\omega_B$ is the dualizing sheaf, and $\shF_B$ is the sheaf of rank two given by the non-split extension
$0\ra \omega_B\ra\shF_B\ra\O_B\ra 0$ corresponding to $H^1(B,\omega_B)=F$.
Similar results for Brauer--Severi curves where obtained by Biswas and Nagaray \cite{Biswas; Nagaraj 2009}
and Novakovi\'c \cite{Novakovic 2012}.

\medskip
The paper is organized as follows:
In Section \ref{Genus-zero curves} and \ref{Quadric hypersurfaces} we collect generalities on genus-zero curves, twisted ribbons, 
and quadric hypersurfaces, which play an important role throughout.
In Section \ref{Surfaces without} we set the stage and examine regular surfaces $X$ without irregularity and bigenus over
arbitrary ground fields $F$, and already formulate our generalization of Iskovskih's Theorem.
The integral curves $C\subset X$ having $h^1(\O_C)=0$ and $C^2\geq 1$ such that all linearly equivalent curves
are  integral are examined in Section \ref{Fields of constants}, and the possible values for $h^0(\O_C)$ and $C^2$ are determined.
This relies on the theory of inseparable pencils, which is established in Section \ref{Inseparable pencils}.
In Section \ref{Bend and break} we start to analyze the crucial new case $C^2=4$,
and formulate the theorem that some linearly equivalent curve $C'$ must be non-integral.
This bend-and-break problem is  quickly reduced to the situation that  $X$ is a modification
of a non-normal quartic surface $V\subset\PP^3$, subject to several arithmetic conditions on   points
with small residue field degree.
In Section \ref{Non-normal quartic} we compute   various  invariants for the normalization maps $Y\ra V$.
It turns out that the branch curve $B\subset V$ is   a space curve of degree three in $\PP^3$,
and the relevant  case is when it is  integral. In Section \ref{Twisted cubics} we then see that $B$ is either
a \emph{twisted cubic}, or is a non-Gorenstein genus-zero curve that arises via denormalization from a projective line
over a cubic field extension $F\subset E$. We cope with such \emph{exotic cubics} via the geometry of the conductor square.
In Section \ref{Blowing-ups} we treat the  twisted cubics $B$ by using the fact that the  blowing-up $\Bl_B(\PP^3)$
unexpectedly carries the structure of a $\PP^1$-bundle over $\PP^2$, and infer that $X$ arises as $\PP^1$-bundle over a genus-zero curve.
To understand the geometry of such bundles, we generalize Grothendieck's Splitting Theorem to
sheaves on regular genus-zero curves in Section  \ref{Locally free}. The final Section \ref{Proofs} completes the proofs of our main results.

\begin{acknowledgement}
This research started when the second  author visited the Universit\'e de Bordeaux, and he wants to thank the host institution 
for its hospitality. 
The first author is supported by the ANR project “FRACASSO” ANR-22-CE40-0009-01.
This research was also conducted in the framework of the   research training group
\emph{GRK 2240: Algebro-Geometric Methods in Algebra, Arithmetic and Topology}, which is funded
by the Deutsche Forschungsgemeinschaft. 
\end{acknowledgement}

\section{Genus-zero curves}
\mylabel{Genus-zero curves}

\newcommand{\codim}{\operatorname{codim}}
\newcommand{\reg}{\text{\rm reg}}
In this section we collect some generalities   that will be used later.   
Throughout, we fix an arbitrary  ground field $F$ of characteristic $p\geq 0$.
The terms \emph{curve} and \emph{surfaces} refer 
to proper schemes that are equi-dimensional, of dimension $n=1$ and $n=2$, respectively.
Each curve $C$ comes with the cohomological invariants
$$
h^0(\O_C)=\dim_F H^0(C,\O_C) \quadand h^1(\O_C)=\dim_F H^1(C,\O_C).
$$
The former is 
the degree of the \emph{ring of constants} $E=H^0(C,\O_C)$, which is a field provided that $C$ is reduced and connected.
In any case,  the integer  $g=h^1(\O_C)$ 
will be  called the \emph{genus}. 
Note that if $\tilde{C}$ is a \emph{twisted form} of some curve $C$, in other words $\tilde{C}\otimes F'\simeq C\otimes F'$ for some
field extension $F\subset F'$, then $h^i(\O_{\tilde{C}})=h^i(\O_C)$.

Following Bayer and Eisenbud \cite{Bayer; Eisenbud 1995}, we say that  a \emph{ribbon} on a curve $C_0$ is a pair $(C,\iota)$, where $\iota:C_0\ra C$ is a closed
embedding having a square-zero sheaf of ideals $\shI$ that is invertible as $\O_{C_0}$-module.
The ribbon  \emph{splits}
if the inclusion $\iota:C_0\ra C$ admits a retraction $\rho:C\ra C_0$; in this case we may regard  $\O_C=\O_{C_0}\oplus\epsilon\shL$
as a \emph{sheaf of dual numbers}, for some invertible $\O_{C_0}$-module $\shL$, and write $C=C_0\oplus\shL$.

Let $f:C'\ra C$ be a birational morphism between curves without embedded components. Then the sheaf of  conductor ideals
$\shI\subset\O_C$, defined as the annihilator for $f_*(\O_{C'})/\O_C$, is also a sheaf of ideals in $\O_{C'}$,
and thus defines both the  branch scheme  $B\subset C$ and the  ramification scheme  $R\subset C'$.
This yields a commutative diagram
\begin{equation}
\label{conductor square}
\begin{CD}
R	@>>>	C'\\
@VVV		@VVf V\\
B	@>>>	C,
\end{CD}
\end{equation} 
which is both cartesian and cocartesian (see  \cite{Fanelli; Schroeer 2020a}, Appendix A). 
In turn, we have a short exact sequence
$0\ra \O_C\ra\O_{C'}\times\O_B\ra \O_R\ra 0$, giving a long exact sequence
\begin{equation}
\label{conductor sequence}
0\ra H^0(\O_C)\ra H^0(\O_{C'})\times H^0(\O_B)\ra H^0(\O_R)\ra H^1(\O_C)\ra H^1(\O_{C'})\ra 0
\end{equation} 
Loosely speaking, the topological space $|C|$ is obtained from $|C'|$ by identifying  the points in $R$ with the same image in $B$,
and the structure sheaf $\O_C$ is the  kernel for $\O_{C'}\times\O_B\ra\O_R$. The 
diagram \eqref{conductor square} is   called \emph{conductor square}, and the passage from $C'$ to $C$ is  called \emph{pinching}.

In what follows, we are interested in the case where the invariants $h^i(\O_C)$ take the minimal possible values,
and find the following locution   useful:

\begin{definition}
\mylabel{genus-zero curve}
A curve $C$ with invariants $h^0(\O_C)=1$ and $h^1(\O_C)=0$ is called a  \emph{genus-zero curve}.
\end{definition}

Note that for each genus-zero curve $C$, there are no embedded components, and the Picard scheme $\Pic_{C/F}$ is \'etale.
Moreover,  each twisted form of $C$ is a genus-zero curve, and 
each   connected modification $C'$ of a reduced subcurve in $C$ is a genus-zero curve, albeit  over the
field extension $F'=H^0(C',\O_{C'})$.

Each finite ring extension $F\subset E$ yields a genus-zero curve via the cocartesian 
\begin{equation}
\label{denormalization of line}
\begin{CD}
\Spec(E) 	@>>>	\PP^1_R\\
@VVV		@VVV\\
\Spec(F)	@>>>	C,
\end{CD}
\end{equation} 
and see from the exact sequence \eqref{conductor sequence} that the pinching $C$ is a genus-zero curve.
A case of particular interest arises when $E$ is a quadratic or cubic  field extension.

To a large extend, the structure of genus-zero curves is determined  by the following result:

\begin{proposition}
\mylabel{classification genus-zero curves}
For each irreducible genus-zero curve $C$ the following holds:
\begin{enumerate}
\item 
We have $C\simeq \PP^1$ if and only if $C$ admits an   invertible sheaf  of degree one.
\item 
If $C$ is regular, the curve is a twisted form of $\PP^1$ or $\PP^1\oplus\O_{\PP^1}(-1)$.
\item 
If $C$ is integral, Gorenstein  and singular, the curve    is obtained from the projective line   over some
quadratic extension $F\subset E$ via the pushout   \eqref{denormalization of line}.
\end{enumerate}
In all three cases, $C$ is a quadric plane curve.  In (ii) the curve must be a  twisted form of $\PP^1$ provided $p\neq 2$.
In   (iii), the curve $C$ is a twisted form
of the union $\PP^1\cup\PP^1$ with $\PP^1\cap \PP^1=\Spec (F)$, or a twisted form of $\PP^1\oplus\O_{\PP^1}(-1)$.
\end{proposition}
 
\proof
Suppose there is an invertible sheaf $\shL$ of degree one. Then $h^0(\shL)\geq \chi(\shL)=\deg(\shL)+\chi(\O_X)=2$.
Choose a global section $s\neq 0$, and let $Z\subset C$ be the zero scheme.
The resulting short exact sequence $0\ra \O_C\ra\shL\ra \shL_C\ra 0$ yields  an exact sequence
$$
H^0(C,\shL)\ra H^0(Z,\shL_Z)\lra H^1(X,\O_X)\lra H^1(X,\shL)\lra 0.
$$
It follows that $\shL$ is globally generated  and $h^1(\shL)=0$, and thus $h^0(\shL)=\chi(\shL)=2$.
The resulting morphism $f:C\ra\PP^1$ has degree one. The  cokernel for $\O_{\PP^1}\ra f_*(\O_C)$ is zero-dimensional
with trivial Euler characteristic, hence $f$ is an isomorphism. This gives (i).

If $C$ is integral and Gorenstein, the dualizing sheaf has $\deg(\omega_C)=-2\chi(\O_C)=-2$.
Arguing as before, we infer that $\shL=\omega_C^{\otimes-1}$ is globally generated with $h^0(\shL)=3$,
and the resulting morphism $f:C\ra\PP^2$ is a closed embedding, with image of degree two.
If $C$ is geometrically reduced, we use (i) after base-changing to $F^\alg$ and see that $C$ is a twisted form
of $\PP^1$. If $C$ is geometrically non-reduced, we regard it as a quadric curve in $\PP^2$,
and see that it is a twisted form of $2L\subset \PP^2$ for some line $L$. This is a ribbon on $\PP^1$
with sheaf of ideals $\shI=\O_{\PP^1}(-1)$. 
From $\Ext^1(\Omega^1_{\PP^1},\shI)=H^1(\PP^1,\O_{\PP^1}(2-1))=0$ we see that the ribbons splits,
so $C$ is a twisted form of $\PP^1\oplus\O_{\PP^1}(-1)$. 
This establishes (ii). 

Suppose that $C$ is integral, Gorenstein and singular. Let $C'\ra C$ be the normalization, and
form the conductor square \eqref{conductor square}. The $C'$ is a regular genus-zero curve over the field extension $F'=H^0(C',\O_{C'})$.
Set
$$
d=[F':F]\quadand l=\dim_{F'} H^0(R,\O_R)\quadand \lambda = dl - \dim_F H^0(B,\O_B).
$$
The short exact sequence \eqref{conductor sequence} yields $1-(d+ (dl-\lambda)) + \lambda$.
The Gorenstein condition ensures $\lambda=dl/2$, see \cite{Fanelli; Schroeer 2020a}, Proposition A.2,
and the equation becomes $d(l-2)=2$. The only solution is $d=2$ and $l=1$, which gives (iii).

The remaining statement arise as follows: From Proposition \ref{quadric singular locus} below  we see that in characteristic $p\neq 2$
are regular quadric curve $C\subset \PP^2$ must be smooth, hence are twisted forms of $\PP^1$.
If $C$ is as in (iii), then the tensor product $E\otimes F^\alg$ is isomorphic to the product
algebra $F^\alg\times F^\alg$ or the ring of dual numbers $F^\alg[\epsilon]$,
producing the union of lines or the ribbon, respectively.
\qed

\begin{corollary}
\mylabel{zariski dense}
Let $C$ be a genus-zero curve that is integral, Gorenstein, and contains at most one rational point.
Then the closed points $c\in C$ of degree two are Zariski dense, and their  local rings $\O_{C,c}$ are regular.
If $C$ is singular, each $\O_C(c)$ generates the Picard group. 
\end{corollary}

\proof
By the proposition, there is a globally generated invertible sheaf $\shL$ of degree two. Choose two global sections
$s_0,s_1$ that generate $\shL$. The resulting morphism $f:C\ra \PP^1$ has degree two.
For each rational point $t\in \PP^1$, the preimage is either  a degree two point, or contains a rational point.
Thus the degree two points $c\in C$ are Zariski dense. 

If $C$ is singular, its description as pinching from $\PP^1_E$ shows that the singular locus
comprises a unique point that is rational, so the $\O_{C,c}$ must be regular.
Since the degree map $\deg:\Pic(C)\ra \ZZ$ is injective and there is no invertible sheaf of degree one,
each $\O_C(c)$ is a generator of the Picard group.
\qed

\medskip
The following terminology  will be useful throughout:
The twisted forms of the projective line are called \emph{Brauer--Severi curves}, and those  having no rational points are 
referred to as  \emph{twisted lines}. The
twisted forms of $\PP^1\cup\PP^1$ with $\PP^1\cap\PP^1=\Spec(F)$ that are integral are called \emph{twisted line pairs}.
If $C$ is a twisted form of $C_0=\PP^1\oplus\O_{\PP^1}(-1)$ that is regular, we say
that $C$ is a \emph{twisted ribbons}; if it is merely reduced and singular we call it
a \emph{twisted ribbon with singularity}. Note that this is only possible over imperfect fields $F$
in character $p=2$, and the former actually needs $\pdeg(F)\geq 2$.

\begin{corollary}
\mylabel{base-change twisted ribbon}
Let $C$ be  a twisted ribbon,  and $F\subset F'$ be an inseparable quadratic extension.
Then the base-change $ C'=C\otimes F'$ is a twisted ribbon, or a twisted ribbon with singularity.
The latter holds   if and only there is a closed point $c\in C$ with $\kappa(c)\simeq F'$.
\end{corollary}

\proof
First note that $p=2$, and that the  base change $C'$ is a twisted form of $\PP^1_{F'}\oplus\O_{\PP^1_{F'}}(-1)$.
It must be integral by \cite{Schroeer 2010}, Lemma 1.3. It regular, it is a twisted ribbon.
If singular, it must be a twisted ribbon with singularity, according to the proposition.  
\qed

\medskip

\begin{lemma}
\mylabel{covering  twisted ribbon}
Let $C$ be a twisted ribbon, and $F\subset E$
be an inseparable quadratic extension. Then there is no  morphisms $f:\PP^1_E\ra C$ of degree two.
\end{lemma}

\proof
Suppose such a morphism exists. Then $f$ is flat, and we get a short exact sequence
$$
0\lra \O_C\lra f_*(\O_{\PP^1_E})\lra \shL\lra 0
$$
for some invertible sheaf $\shL$ on the twisted ribbon $C$. From $\chi(\O_C)=1$ and $\chi(\O_{\PP^1_E})=2$ we get
$\deg(\shL)=\chi(\shL)-\chi(\O_C)=0$, and thus $\shL=\O_C$. Using $\Ext^1(\O_C,\O_C)=H^1(C,\O_C)=0$,
we see that the above short exact sequence of $\O_C$-modules splits. In turn, the ring structure takes the form 
 $f_*(\O_{\PP^1_E})=\O_C[T]/(P)$ for some polynomial $P(T)=T^2+\alpha T+\beta$ whose  coefficients belong to 
$H^0(C,\O_C)=k$. It follows that $\PP^1_E=C\otimes_FE'$ for some quadratic field extension $F\subset E'$.
Taking global sections we obtain $E'=E$.  By Proposition \ref{base-change twisted ribbon}, the base change  $C\otimes E$ is a twisted ribbon, perhaps with singularity,
By Proposition \ref{classification genus-zero curves} it contains at most one rational point, contradiction.
\qed
 
\begin{lemma}
\mylabel{covering implies line}
Let $f:C'\ra C$ be a  finite flat morphism of degree two between genus-zero curves.
If $C$ is integral, it must be isomorphic to $\PP^1$.
\end{lemma}

\proof
By flatness, the cokernel $\shL=f_*(\O_C')/\O_C$ is invertible, with Euler characteristic $\chi(\shL)=\chi(\O_{C'})-\chi(\O_C)=0$.
Then $\deg(\shL)=\chi(\shL)-\chi(\O_C)=-1$, and the assertion follows from Proposition \ref{classification genus-zero curves}.
\qed

\section{Quadric hypersurfaces}
\mylabel{Quadric hypersurfaces}

We keep the set-up from the previous section, and 
now examine \emph{quadric hypersurfaces} $V\subset\PP^{n+1}$ of arbitrary dimension $n\geq 0$.
They  are  defined by a non-zero homogeneous  polynomial  $\Phi$ of degree two in $n+2$ indeterminates  with coefficients from $F$.
After a change of coordinates, we may write
\begin{equation}
\label{quadric polynomial}
\Phi  = \sum_{i=0}^m \lambda_iT_i^2
\end{equation} 
in the indeterminates $T_i$, for some $0\leq m\leq n+1$ and     $\lambda_i\in F^\times$, provided $p\neq 2$.
In characteristic two, we instead may write
\begin{equation}
\label{quadric polynomial p=2}
\Phi = \sum_{i=1}^r (X_iY_i + \alpha_iX_i^2 + \beta_iY_i^2) + \sum_{j=1}^s \gamma_jZ_j^2
\end{equation} 
in some indeterminates  $X_i,Y_i,Z_j$, with $0\leq 2r+s-1\leq n+1$ and    $\alpha_i,\beta_i,\gamma_j\in F$, where the 
$\gamma_1,\ldots,\gamma_s$ are  linearly independent over the subfield $F^p$, and $r=0$ or $s=0$ means
that the respective sums disappear. See \cite{Arf 1941}, Satz 2  for more details.

Note that our designations $T_i$ or $X_i,Y_i,Z_j$ for the indeterminates will indicates that the characteristic is $p\neq 2$ or $p=2$, 
respectively. In the latter case, let $F^p\subset E$ be the height-one extension generated by the fractions $\gamma_j/\gamma_k$, $1\leq j,k\leq s$,
and consider the ensuing number 
$$ 
\delta_\text{\rm reg}= \pdeg(E/F^p).
$$
Recall that the \emph{$p$-degree} is the cardinality of a $p$-basis of $E$ over $F^p$, or equivalently of  a basis for
$\Omega^1_{E/F^p}$ over $E$. 
The  \emph{jacobian ideal} $\mathfrak{j}=(X_1,Y_1,\ldots,X_r,Y_r)$ defines a linear subspace $\PP^{s-1}$;
for  $s\geq 1$ we write  $V_0\subset \PP^{s-1}$ for the induced quadric hypersurface, which is defined by the 
the Fermat polynomial   $\Phi_0= \sum_{j=1}^s \gamma_jZ_j^2$.   
 
\begin{proposition}
\mylabel{quadric singular locus}
In the above situation, the 
scheme of non-smoothness $\Sing(V/F)$ is the   subscheme in $\PP^{n+1}$ defined by  the respective homogeneous ideals 
$$
\ideala=(T_0,\ldots,T_m)\quadand \idealb=(X_1,Y_1,\ldots,X_r,Y_r, \gamma_1Z_1^2 + \ldots + \gamma_sZ_s^2).
$$
This coincides with the singular locus $\Sing(V)$, except for $p=2$ and $s\geq 1$. In this case,  
 $\Sing(V)\subset V$ has codimension $2r+\delta_\reg$, and coincides with the singular locus of the 
Fermat hypersurface $V_0\subset\PP^{s-1}$.
\end{proposition}

\proof
Computing partial derivatives, one immediately gets the statement on the scheme of non-smoothness.
The second statement is a special case of \cite{Schroeer 2010}, Theorem 3.3. 
\qed

\medskip
For $p\neq 2$ we thus see that $\Sing(V/F)=\Sing(V)$. If non-empty, this locus has codimension $0\leq m\leq n$ 
inside the $n$-dimensional quadric hypersurface $V$.  In characteristic two, we get
$$
\codim_V\Sing(V/F) = 2r\quadand 
\codim_V\Sing(V)= 2r+\delta_\reg.
$$
Using Serre's Criterion for normality (\cite{EGA IVb}, Theorem 5.8.6), we immediately get:  

\begin{corollary}
\mylabel{quadric   normal}
The quadric hypersurface $V\subset\PP^{n+1}$ is reduced if and only if $m\geq 1$  or   $r\geq 1$ or $\delta_\reg\geq 1$.
It is normal if and only if $m\geq 2$ or $r\geq 1$ or $\delta_\reg\geq 2$.
\end{corollary}

The following observations are   useful throughout:

\begin{lemma}
\mylabel{connectedness and normalization}
The singular locus $\Sing(V)$ is connected.
If $V$ is reduced and non-normal, the normalization is isomorphic to $\PP^n_E$,
where $F\subset E$ is either a quadratic field extension or the product ring $E=F\times F$. 
\end{lemma}

\proof
The connectedness statement immediately follows from Proposition \ref{quadric singular locus}, except in the case $p=2$ and $s\geq 1$.
It then suffices to treat the case    
$$
\Phi=\sum_{j=1}^{s}\gamma_jZ_j^2\quadand s=n+2\quadand r=0.
$$
Set $t=\delta_\reg$. 
Without loss of generality we may assume that $\gamma_1,\ldots,\gamma_t\in F$ are $p$-linearly independent over $F^p$,
and that $\gamma_{t+1}=1$.
By \cite{A 4-7},  Chapter V, \S 13, No.\ 2, Theorem 1 
there are  derivations $D_i:F\ra F$ with $D_i(\gamma_j)=\delta_{ij}$  for $1\leq i,j\leq t$,
and thus
$$
D_i(\Phi)= Z_i^2+\sum_{j=t+2}^s D_i(\gamma_j) Z_j^2.
$$
Set $\idealc=(\Phi,D_1(\Phi),\ldots,D_t(\Phi))$. The resulting closed subscheme $V_+(\idealc)\subset\PP^{n+1}$
is isomorphic to the   hypersurface in  $\Proj F[Z_{t+1},\ldots,Z_{n+2}]$
defined by a Fermat equation of the form $Z_{t+1}^2+\sum_{j=t+2}^s\gamma'_jZ_j^2=0$, hence
is irreducible of dimension $(n+2)-(t+1)= n+1-\delta_\reg$. It contains $\Sing(V)$, which by Proposition \ref{quadric singular locus} 
has the same dimension. Thus the inclusion $\Sing(V)\subset V_+(\idealc)$ is an equality,
and   connectedness follows.

Now suppose that $V$ is reduced and non-normal, with normalization $V'$.
Assume first  $p\neq 2$. By Corollary \ref{quadric normal}, it suffices to treat the case $\Phi=T_0^2-\lambda T_1^2$
for some $\lambda\in F^\times$, which in characteristic two   does not belong to $F^{2\times}$.
Then the ring $E=F[U]/(U^2-\lambda)$ is regular, and we write $\omega\in E$ for the class of the indeterminate $U$. 
As for \cite{Schroeer 2010}, Proposition 4.1, one easily checks
that the homomorphism
$$
F[T_0,\ldots,T_{n+1}]\lra E[T'_1,\ldots,T'_{n+1}]
$$
of graded $F$-algebras  given by $T_0\mapsto \omega T'_1$ and $ T_i\mapsto T'_i$ for $1\leq i\leq n+1$
induces a birational morphism $\PP^n_E\ra V\subset\PP^{n+1}$, hence $V'=\PP^n_E$.
\qed

\medskip

\section{Surfaces without irregularity and bigenus}
\mylabel{Surfaces without}

Let $F$ be a ground field of characteristic $p\geq 0$, not necessarily perfect,
 and $X$ be a     regular   surface, not necessarily smooth, with numerical invariants
$$
h^0(\O_X)=1\quadand h^1(\O_X)=  h^0(\omega_X^{\otimes 2}) =0.
$$
Clearly, we have $h^0(\omega_X)=0$, and Serre Duality ensures $h^2(\O_X)=0$, such that $\chi(\O_X)=1$.
The Chern number  $c_1^2=(\omega_X\cdot\omega_X)$ and the Picard number $\rho\geq 1$ are further invariants
of interest. The same goes for   invariants stemming from the coherent sheaf $\Omega^1_{X/F}$, which is locally free of rank two
if and only if $X$ is smooth. 

Throughout we also assume that $X$ is \emph{minimal}. In other words, there is no integral curve $R\subset X$ 
such that 
$$
(R\cdot R)=(R\cdot\omega_X)=-h^0(\O_E).
$$
Such $R$ are isomorphic to the projective line over $E=H^0(R,\O_R)$, 
and are called \emph{exceptional curves of the first kind}; one may   call them \emph{$(-1)$-curves} as well, in reference to the sign occurring the above equations.
We can already state the main result of this paper:

\begin{theorem}
\mylabel{trichotomy}
In the above situation, one of the following conditions holds:
\begin{enumerate}
\item The surface $X$ is isomorphic to  $\PP^2$,   or  to a   quadric    in $\PP^3$.
\item There is a morphism $f:X\ra B$ with $f_*(\O_X)=\O_B$ such that the base  and the generic fiber   are regular genus-zero curves.
\item The dualizing sheaf $\omega_X$ generates the Picard group $\Pic(X)$.
\end{enumerate}
\end{theorem}

This trichotomy generalizes Iskovskih's result (\cite{Iskovskih 1980}, Theorem 1), 
by replacing smoothness with regularity and allowing imperfect base fields. 
The proof requires extensive preparation, and will be completed in Section \ref{Proofs}.
Note that our surface $X$ may be geometrically non-normal, or even geometrically non-reduced.
In this general setting, one has  to cope with two new issues: 
The integral curves $C\subset X$ may have larger self-intersection $C^2$, compared to the smooth case,
and the field of constants $E=H^0(C,\O_C)$ could be an extension of the ground field $F$.

In this section we closely follow  the exposition of Koll\'ar, Smith and Corti (\cite{Kollar; Smith; Corti 2004},
Section 3.2), but with special attention to the above issues.
Let us start with  the following key fact, classically called ``Termination of Adjunction''.

\begin{lemma}
\mylabel{termination adjunction}
For each invertible sheaf $\shL$ we have $h^0(\shL\otimes\omega_X^{\otimes n})=0$
for all sufficiently large $n\geq 0$.
\end{lemma}

\proof
Serre Duality ensures $h^2(\omega_X^{\otimes-1})= h^0(\omega_X^{\otimes 2})$, which vanishes by our assumption.
Thus Riemann--Roch gives
$h^0(\omega_X^{\otimes-1})\geq \chi(\omega_X^{\otimes-1}) = c_1^2 +\chi(\O_X) = c_1^2 + 1$.
We have $\omega_X^{\otimes-1}\not\simeq\O_X$ because $h^0(\omega_X)=0$.
So if 
$c_1^2\geq 0$, then   each ample divisor $H\subset X$ has $(\omega_X\cdot H)<0$,
therefore $(\shL\otimes\omega_X^{\otimes n}\cdot H)<0$ for all $n$ sufficiently large, and then $h^0(\shL\otimes\omega_X^{\otimes n})=0$.

Suppose now $c_1^2<0$. Then $(\shL\otimes\omega^{\otimes n}_X\cdot\omega_X)<0$  for all  $n\geq n'$ with some integer $n'$, and then the 
$\shL\otimes\omega^{\otimes n}_X$ are non-trivial.
Seeking a contradiction, we assume that
for  infinitely many  $n_i\geq 0$, $i\geq0 $ the invertible sheaf $\shL\otimes\omega_X^{\otimes n_i}$ admits a non-zero  global section $s_i$.
Without loss of generality we may assume $n_i\geq n'$.
Let $D_i\subset X$ be the zero scheme of the $s_i$. Decompose $D_0=\sum m_jC_j$ into irreducible components.
Then some $C=C_j$ has $(C\cdot \omega_X)<0$.  Since $X$ is minimal, we must have $C^2\geq0$,
and thus $(\shL\otimes\omega_X^{\otimes n_i}\cdot C)\geq 0$ for all $i\geq 0$. On the other hand, $(\omega_X\cdot C)<0$
ensures that for sufficiently large $i\geq 0$ we have  $(\shL\otimes\omega_X^{\otimes n_i}\cdot C)<0$, contradiction.
\qed

\medskip
This has an immediate consequence for the Picard scheme:

\begin{proposition}
\mylabel{picard scheme}
The canonical map $\Pic_{X/k}\ra\Num_{X/k}$ of group schemes is an isomorphism.
\end{proposition}

\proof
The Lie algebra for $\Pic_{X/k}$ is the cohomology group $H^1(X,\O_X)$. This vanishes by assumption,
so $\Pic^0_{X/k}$ is trivial. To see that the local system $\Num_{X/k}$ contains no torsion,
it suffices to treat the case that $k$ is separably closed. Seeking a contradiction, we assume
that there is an invertible sheaf $\shL$ of finite order $m>1$ in the Picard group.
Then $\shL$ is numerically trivial with $h^0(\shL)=0$, so Serre Duality and Riemann--Roch gives
$h^2(\shL)\geq \chi(\shL)=\chi(\O_X)=1$.
Again by Serre duality, the  invertible sheaf $\shN=\omega_X\otimes\shL^\vee$ admits a non-zero global section,
and the same holds for the tensor power $\shN^{\otimes m}=\omega_X^{\otimes m}$.
Thus $h^0(\omega_X^{\otimes n})\neq 0$ for all positive multiples $n$ of $m$, in contradiction to Lemma \ref{termination adjunction}.
\qed

\medskip
We conclude that $\Pic(X)=\Num(X)$ is a free abelian group. Furthermore, the class of the dualizing
sheaf $\omega_X$ is non-zero, in light of  Lemma \ref{termination adjunction}.

\begin{proposition}
\mylabel{good curve}
Suppose the dualizing sheaf $\omega_X$ does not generate the Picard group $\Pic(X)$.
Then there is an integral curve $C\subset X$ such that 
$h^1(\O_C)=0$ and $ C^2\geq 0$,
and that every linearly equivalent   $C'\subset X$ is also integral.
\end{proposition}

\proof
For each ample divisor  $H\subset X$, we   consider  the function     $n\mapsto h^0(\omega_X^{\otimes n}(H))$.
It is  non-zero for $n=0$, but vanishes for $n$ sufficiently large, by Lemma \ref{termination adjunction}.
Let $n=n_H$ be the largest integer  such that $\shN=\omega_X^{\otimes n}(H)$ has a non-zero global section 
but $\shN\otimes\omega_X=\omega_X^{\otimes n+1}(H)$ does not. 
Fix such a global section $s$, and regard it as homomorphism $s:\omega_X^{\otimes -n}\ra\O_X(H)$.
Now recall that $\Pic(X)$ is generated by the classes of very ample divisors.
If the $s$ are bijective for all $H$, then $\omega_X$ would generate the Picard group, contradiction.
Therefore for some $H$ the map $s$ is not bijective.

Summing up, we find some non-trivial invertible sheaf $\shN$ with $h^0(\shN)\neq 0$ but $h^0(\shN\otimes\omega_X)=0$.
Now re-choose a non-zero $s\in H^0(X,\shN)$ so that the zero scheme $D\subset X$ maximizes $\sum_{i=1}^rm_i$,
where $D=\sum_{i=1}^r m_iC_i$ is the decomposition into irreducible components.
Let $C=\sum m_i'C_i$ be any subcurve of $D$.
The resulting invertible sheaf $\shL=\O_X(C)$
has the property that $h^0(\shL^{\otimes-1})=0$ and $h^0(\shL\otimes\omega_X)=0$.
Now Serre Duality ensures $h^2(\shL^{\otimes-1})=0$, and Riemann--Roch gives
\begin{equation}
\label{estimate}
0=h^0(\shL^{\otimes-1}) + h^2(\shL^{\otimes-1}) \geq \chi(\shL^{\otimes-1}) = \frac{(C+K_X)\cdot C}{2} +\chi(\O_X).
\end{equation}
The Adjunction Formula yields $\deg(\omega_C)=(C+K_X)\cdot C \leq -2\chi(\O_X)=-2$, whereas Serre Duality
gives $\deg(\omega_C)=-2\chi(\O_C)$, thus $\chi(\O_C)\geq 1$. 

Now additionally assume that for the subcurve $C\subset D$,  the finite $k$-algebra  $E=H^0(C,\O_C)$ is a field. The latter
condition holds, for example, if $C$ is reduced and connected, and in particular if $C$ is integral.
The cohomology group  $H^1(C,\O_C)$ carries the structure of an $E$-vector space. We have 
$$
1\leq \chi(\O_C) = [E:k](1 - \dim_E H^1(C,\O_C))
$$
and conclude $h^1(\O_C)=0$. In particular, $C$ is a genus-zero curve over the extension field $E$. 

The preceding paragraph applies in particular to  $C=C_i$,  so each $C_i$ is a genus-zero curve over the field $E_i=H^0(C_i,\O_{C_i})$.
Seeking a contradiction, we assume now that $C_i^2<0$ for all $1\leq i\leq r$. 
If $C_i^2=-[E_i:k]$ then $C_i\subset X$ would be an exceptional curve of the first kind, contradiction.
Thus $C_i^2\leq -2[E_i:k]$, and from 
$$
-2[E_i:k]=\chi(\O_{C_i}) = \deg(\omega_{C_i}) = C_i^2+(\omega_X\cdot C)\leq -2 [E_i:k]  + (\omega_X\cdot C),  
$$
we infer $(\omega_X\cdot C_i)\geq 0$, and thus $(\omega_X\cdot D)\geq 0$. 
Now recall that $\omega_X^{\otimes -n}(D)=\O_X(H)$ is ample, and therefore
$$
(\omega_X(D)\cdot D) = (\omega_X^{\otimes-n}(D)\cdot D)  + (\omega_X^{\otimes n+1}\cdot D) \geq (\omega_X^{\otimes n+1}\cdot D)\geq 0.
$$
On the other hand, we saw in \eqref{estimate} with $C=D$ that $(\omega_X(D)\cdot D)<0$, contradiction.
This shows that some $C=C_j$ has $h^1(\O_C)=0$ and $C^2\geq 0$. By maximality of $\sum_{i=1}^rm_i\geq 1$,
every linearly equivalent curve $C'$ remains integral.
\qed

\begin{proposition} 
\mylabel{globally generated}
Let  $C\subset X$ be an integral curve with $h^1(\O_C)=0$ and $C^2\geq 0$.
Then   the invertible sheaf $\shL=\O_X(C)$ is globally generated, with $h^0(\shL)= 1+h^0(\O_C) + C^2$.
\end{proposition}

\proof
The short exact sequence $0\ra \O_X\ra\shL\ra\shL_C\ra 0$ gives a long exact sequence
$$
0\lra H^0(X,\O_X)\lra H^0(X,\shL)\lra H^0(C,\shL_C)\lra H^1(X,\O_X),
$$
and the term on the right vanishes. So to see that $\shL$ is globally generated, it suffices to verify
this for $\shL_C$. This indeed holds, because $C$ is an integral genus-zero curve over the
field $E=H^0(C,\O_C)$, and $\deg(\shL_C)=C^2\geq 0$. Furthermore, we have $h^1(\shL_C)=0$, and thus
$h^0(\shL_C)=\chi(\shL_C) = \deg(\shL_C)+\chi(\O_C) = C^2+h^0(\O_C)$.
The assertion on $h^0(\shL)$   follows from the above long exact sequence.
\qed

\medskip
In the above situation, the  globally generated invertible  sheaf $\shL=\O_X(C)$ defines a morphism 
$$
f:X\lra\PP^n\qquad\text{with}\qquad f^*\O_{\PP^n}(1)=\shL,
$$
where   $n=h^0(\O_C)+C^2$. The schematic image $V\subset\PP^n$ is an integral closed subscheme of dimension one or two.
This image is not a linear subscheme, because the induced map $H^0(\PP^n,\O_{\PP^n}(1))\ra H^0(X,\shL)$
is bijective. Note that the image    may or may not be regular, and that we have the formula 
\begin{equation}
\label{degree equation}
C^2=\deg(f)\cdot\deg(V),
\end{equation} 
where the degree for $f:X\ra V$ is formed in the sense of Kleiman (\cite{Kleiman 1966}, Chapter I, Section 2, Definition on page 299).
The following already settles parts of  Theorem \ref{trichotomy}:

\begin{proposition}
\mylabel{small self-intersection}
Suppose there is an integral curve $C\subset X$ with  $h^1(\O_C)=0$ and   $0\leq C^2\leq 2$. Then one of the following holds:
\begin{enumerate}
\item The surface $X$ is isomorphic to the projective plane $\PP^2$, or isomorphic to a   quadric surface in $\PP^3$.
\item There is morphism $g:X\ra B$ with $\O_B=g_*(\O_X)$ such that the base $B$ and and the generic fiber $X_\eta$ are regular genus-zero curves.
\end{enumerate}
Moreover, $X$ is geometrically normal provided the characteristic is  $p\neq 2$.
\end{proposition}

\proof
Suppose first that $C^2=0$. Then the image of $f:X\ra \PP^n$ must be a curve. Let $B=\Spec f_*(\O_X)$
be the Stein factorization, and $g:X\ra B$ be the resulting morphism. Then $B$ is integral and normal.
The Leray--Serre spectral sequence gives $H^0(B,\O_B)=H^0(X,\O_X)$ and an exact sequence
$$
0\lra H^1(B,\O_B)\lra H^1(X,\O_X)\lra H^0(B,R^1g_*(\O_X))\lra 0.
$$
Consequently  $B$ is a regular genus-zero curve and $R^1g_*(\O_X)$ is locally free. 
Since $C^2=0$ the image $f(C)\subset\PP^n$ is a singleton,
and thus $C=g^{-1}(b)$ for some closed point $b\in B$. 
The formation of $R^1g_*(\O_X)$ commutes with base-change, because $g$ is flat with one-dimensional fibers,
and thus  $R^1g_*(\O_X)\otimes \kappa(b)=H^1(C,\O_C)=0$. Semicontinuity now gives $R^1g_*(\O_X)\otimes \kappa(\eta)=0$,
in other words, the generic fiber $X_\eta$ is a genus-zero curve.

For the remaining cases the image $V\subset\PP^n$ is an integral surface, 
subject to \eqref{degree equation}. Also note  that the  integer $h^0(\O_C)$ divides   $ (\shL\cdot C)=C^2$.
Suppose now that $C^2=1$. Then  $f:X\ra V$ is birational, $h^0(\O_D)=1$ and thus $h^0(\shL)=3$.
It follows   $V=\PP^2$. This gives a birational morphism $f:X\ra\PP^2$, which is a sequence of blowing-ups.
It must be an isomorphism because $X$ is minimal.

It remains to treat the case $C^2=2$. We then have $1\leq h^0(\O_C)\leq 2$ and thus $4\leq h^0(\shL)\leq 5$.
To cope with the possible constant field extension, we observe that $\shL|C$ is generated by two global sections.
We thus find $s_0,\ldots,s_3\in H^0(X,\shL)$ without common zero, and now replace   $f$
by the morphism $f:X\ra\PP^3$ with $f^{-1}(t_i)\otimes 1=s_i$, where $\PP^3=\Proj F[t_0,\ldots,t_3]$.
Now the image $V$ must be a quadric surface, and the induced  $f:X\ra V$ is birational.
By the argument in the preceding paragraph, $f$ is an isomorphism provided that $V$ is regular. 

Seeking a contradiction, we assume that the quadric surface $V$ is not regular. The morphism $f:X\ra V$ factors over
the normalization $V'\ra V$, thus $h^0(\O_{V'})=1$. It follows from Lemma \ref{connectedness and normalization} that $V$ is normal,
and the singular locus is a singleton $v_0\in V$. It follows from  the description in Proposition \ref{quadric singular locus}
that this point is rational.
Choose some global section of $\O_V(1)$ so that the zero scheme $D\subset V$ passes through $v_0$.
Then the preimage $f^{-1}(D)$ is linearly equivalent to $C$ but reducible, contradiction.

It remains to check that $X$ is geometrically normal, provided $p\neq 2$.
If $C^2=0$, then the base and the generic fiber of the fibration $g:X\ra B$
are smooth, according to Proposition \ref{classification genus-zero curves}, so the   $\Sing(X/F)$ is finite.
By Serre's Criterion (\cite{EGA IVb}, Theorem 5.8.6) the scheme $X$ is geometrically normal.
There is nothing to prove if $X=\PP^2$, so it remains to treat the case that $X\subset\PP^3$ is a regular quadric.
Then $X$ is smooth, by Proposition \ref{quadric singular locus}.
\qed

\section{Fields of constants and  self-intersections}
\mylabel{Fields of constants}

We keep the set-up from the previous section, such that $X$ is a  minimal
regular surface over a ground field $F$ of characteristic $p\geq 0$,   with $h^0(\O_X)=1$ and $h^1(\O_X)=h^0(\omega_X^{\otimes 2})=0$. 
We also assume that there is an integral curve $C\subset X$
with $h^1(\O_C)=0$ and self-intersection  $C^2>0$, such that all linearly equivalent curves $C'$  are  also integral.
By Proposition \ref{globally generated} the invertible sheaf $\shL=\O_X(C)$ is globally generated and has  
\begin{equation}
\label{dimension general system}
h^0(\shL)=1+h^0(\O_C)+C^2\geq 3.
\end{equation} 
In this section, we will determine the second  summand,    and constrain the third.  

\begin{proposition}
\mylabel{field of constants}
The degree of the constant field extension is   $h^0(\O_C)=1$.
\end{proposition}

\proof 
Seeking a contradiction, we assume that $E=H^0(C,\O_C)$ is a non-trivial extension of the ground field $F$.
Using $C^2>0$ and the Hodge Index Theorem  (\cite{Grothendieck 1958}, Theorem 1.1), we see that $C$ is geometrically connected,
thus $F\subset E$ is purely inseparable (\cite{A 4-7}, Chapter V, \S7, No.\ 8, Proposition 13). 
Hence we are in characteristic $p>0$ and the field $F$ is imperfect.
By   Theorem \ref{pencils characteristic two}  below we have  $p=2$ and $[E:F]=2$.  

Regarding $C$ as a genus-zero curve over $E$,
we find with Corollary \ref{zariski dense} a closed point $z\in C$ such that the residue field $K=\kappa(z)$ has degree
at most two over $E$, and that the local ring $\O_{C,z}$ is regular.
Choose a global section of $\shL|C$  that vanishes at $z\in C$ and has only simple zeros, and extend it to a global section $s'\in H^0(X,\shL)$.
Let $C'\subset X$ be its zero-scheme and set $Z=C\cap C'$. The ensuing blowing-up $r:\Bl_Z(X)\ra X$ comes
with a fibration $h:\Bl_Z(X)\ra\PP^1$, as explained in the next section. Again by Theorem \ref{pencils characteristic two}, 
the Stein factorization $B$ is a twisted line or twisted ribbon. In particular, our ground field $F$ is infinite. Moreover, the projection $g:B\ra \PP^1$ is radical of  degree two,
and for each rational point $t\in\PP^1$, the fiber $g^{-1}(t)$ is reduced.

Suppose   that $B$ is a twisted line.  
First note that over some dense open set  the formation of $h_*(\O_X)$ commutes with arbitrary base-change.
By Corollary \ref{zariski dense} we can pick a closed point $b\in B$ that is contained  in this open set and whose residue field $L=\kappa(b)$ is a separable quadratic extension.
Then  the image $g(b)\in \PP^1$ is not rational, hence also has residue field $L$.
Consider the corresponding $L$-valued points $b\in B\otimes_FL$ and $g(b)\in\PP^1_L$. The base-change $B\otimes_F L$ is the projective line over $L$,
and the fiber over $t=g(b)$ has coordinate ring $L[\epsilon]$.
So this point corresponds to a divisor $D$ on the base-change $X\otimes L$ that    is linearly equivalent to $C\otimes L$
and has $H^0(D,\O_D)=L[\epsilon]$, and the structure morphism $D\ra\Spec L[\epsilon]$ is flat.
Since $X\otimes L$ remains regular, the subcurve $H\subset D$ defined by $\epsilon=0$ is Cartier, with $D=2H$, and thus $H^2>0$.
 On the other hand, the  short exact sequence $0\ra \epsilon\O_H\ra \O_D\ra \O_H\ra 0$ 
gives $H^2=-\deg(\epsilon\O_H)=0$, contradiction.
 
It remains to treat the case that $B$ is a twisted ribbon.
The fiber $r^{-1}(z)$ of the blowing-up is a copy of $\PP^1_K$, and the degree of  composite map
$$
\PP^1_K\lra B\lra \PP^1
$$
is the degree of the finite scheme $r^{-1}(z)\cap h^{-1}(0)=\Spec(K)$. By construction  $[K:F]=[K:E]\cdot [E:F]\leq 4$.
Using Lemma \ref{covering  twisted ribbon} we see $\deg(\PP^1_K/B)>2$, which gives $[K:F]>4$, contradiction. 
\qed

\medskip
So   formula \eqref{dimension general system}   simplifies to $h^0(\shL)=2+C^2$. It turns out that there are only three possibilities: 

\begin{proposition}
\mylabel{selfintersection numbers}
The self-intersection number satisfies  $C^2\mid 4$. If there are two rational points $a\neq b$ on $X$ we actually have $C^2\mid 2$.
\end{proposition}

\proof
First suppose that there are rational points $a\neq b$, and consider the finite subscheme $Z\subset X$
with coordinate ring $\O_{X,a}/\maxid_a\times \O_{X,b}/\maxid_b^2$, which has $h^0(\O_Z)=1+3=4$.
If $C^2\geq 3$ we have $h^0(\shL)\geq 5$, and find some curve $C'$ that is linearly equivalent to $C$ and contains $Z$.
Then $C'$ is a genus-zero curve  that is integral, Gorenstein and singular with two rational points, 
in contradiction to Proposition \ref{classification genus-zero curves}.
This settles the second assertion.

Suppose now that $X$ contains at most one rational point. By Corollary \ref{zariski dense}, the integer $C^2\geq 1$ is even.
Seeking a contradiction, we now assume that $C^2\geq 6$, hence $h^0(\shL)\geq 8$.
Fix a closed point $c\in C$ of degree two and form the finite subscheme $Z\subset X$ with coordinate ring $\O_{X,c}/\maxid_c^2$,
which has $h^0(Z)=6$. Again we find a curve $C'$ that is linearly equivalent to $C$ and contains $Z$.
By construction, the local ring $\O_{C',c}$ is singular and the residue field $\kappa(c)$ has degree two,
contradiction.
\qed

\section{Inseparable pencils}
\mylabel{Inseparable pencils}

In this section we establish some general facts, which appear to be of independent interest  and 
have been used in the previous section.
Let  $F$ be a  ground field of characteristic $p\geq 0$, and $X$ be a proper normal scheme
with $h^0(\O_X)=1$, of dimension $n\geq 2$.  Let $\shL$ be an invertible sheaf not isomorphic to $\O_X$.
For each non-zero global section $s\in H^0(X,\shL)$, the resulting short exact sequence
\begin{equation}
\label{sequence defining divisor}
0\lra \shL^{\otimes-1}\stackrel{s}{\lra}\O_X\lra \O_D\lra 0
\end{equation} 
defines an effective Cartier divisor $D\subset X$.  The goal of this section is to understand the finite $F$-algebra 
$E=H^0(D,\O_D)$ in dependence of the global section $s\neq 0$. We start with the following observation:

\begin{lemma}
\mylabel{independent degrees}
Suppose  $h^1(\O_X)=0$. Then  the number  $h^0(\O_D)=[E:F]$ does not depend on the global section $s\neq 0$.
\end{lemma}

\proof
The short exact sequence \eqref{sequence defining divisor}
yields a long exact sequence
$$
H^0(X,\shL^{\otimes-1})\lra H^0(X,\O_X)\lra H^0(D,\O_D)\lra H^1(X,\shL^{\otimes-1})\lra H^1(X,\O_X).
$$
The term on the right vanishes by assumption, and the term on the left is zero
because $X$ is integral and $\shL\not\simeq \O_X$ admits a non-zero global section.
Thus   $h^0(\O_D)=1+h^1(\shL^{\otimes-1})$, which obviously   does not depend on the global section.
\qed

\medskip
Suppose now we have two non-zero global sections  $s_0,s_1\in H^0(X,\shL)$ such that the resulting common zero-scheme $Z\subset X$
has  codimension two, and is also reduced.
Write $D_i\subset X$ with  $i=0,1$ for the resulting effective Cartier divisors.
Each irreducible component of $D_i$ has codimension one, according to Krull's Principal Ideal Theorem.
Moreover, there are no common irreducible components, because $Z=D_1\cap D_2$ has codimension two.
Being normal, the scheme $X$ satisfies Serre's Condition $(S_2)$, thus $D_i$ satisfies $(S_1)$,
and we conclude that $Z$ is an effective Cartier divisor inside $ D_i$.
Thus at each point $z\in Z$, the sheaf of ideals $\shI\subset\O_X$ for the closed subscheme $Z$
is generated by a regular sequence contained in the maximal ideal $\maxid_z\subset\O_{X,z}$.
For each generic point $\zeta\in Z$, the local ring $\O_{Z,\zeta}$ and hence also $\O_{X,\zeta}$ are regular.

Let $V=\Bl_Z(X)$ be the blowing-up,  and write $r:V\ra X$ for the resulting morphism.
First note that 
$r_*(\O_V)=\O_X$ and $R^ir_*(\O_V)=0$ for $i\geq 1$, according  to \cite{SGA 6}, Expos\'e VII, Lemma 3.5,
and thus  
$h^j(\O_V)=h^j(\O_X)$ for all $j\geq 0$.
According to \cite{Lorenzini; Schroeer 2023}, Proposition 3.4, the scheme $V$ satisfies Serre's Condition $(S_2)$, 
and furthermore is regular over an open set
$U\subset X$ that contains all codimension one-points in $X$ and all generic points in $Z$. It follows that for  each point $\xi\in V$
of codimension one, the local ring $\O_{V,\xi}$ is regular, hence $V$ is normal. Moreover, the 
exceptional divisor $R=g^{-1}(Z)$ is reduced. 

Consider the invertible sheaf $\shN=r^*(\shL)(1)=r^*(\shL)(-R)$.
The   short exact sequence $0\ra \shN\ra r^*(\shL)\ra r^*(\shL)|R\ra 0$ yields an exact sequence
$$
0\lra H^0(V,\shN)\lra H^0(V,r^*(\shL))\lra H^0(R,r^*(\shL)|R).
$$
In turn, the  $s\in H^0(X,\shL)$ that vanishes along $Z$ can also be seen as the elements in $H^0(V,\shN)$.
Each such $s\neq 0$ thus defines   effective Cartier divisors $D\subset X$ and $D'\subset V$, 
where the latter is the strict transform of the former,
and the induced map $r:D'\ra D$ is the blowing-up with center $Z\subset D$.
Since $Z$ is already Cartier   inside $D$, this gives an identification $D'=D$.

The strict transforms $D'_0$ and $D'_1$ are disjoint, which  follows from \cite{Perling; Schroeer 2017}, Lemma 4.4
and our assumption that $Z$ is reduced. 
In turn, the invertible sheaf $\shN$ is globally generated by $s_0$ and $s_1$, which define 
a morphism 
$$
h:V\lra\PP^1\quad\text{with}\quad h^*(\O_{\PP^1}(1))=\shN.
$$
The morphism is surjective, because  $\shN$ is numerically non-trivial, and therefore flat.
Let  $B=\Spec h_*(\O_V)$ be the Stein factorization, with resulting morphisms $f:V\ra B$ and $g:B\ra\PP^1$.
Then $\O_B=f_*(\O_V)$, so $B$ is normal, with $h^0(\O_B)=1$. Moreover, the Leray--Serre spectral sequence
gives an exact sequence
\begin{equation}
\label{five term exact sequence}
0\lra H^1(B,\O_B)\lra H^1(V,\O_V)\lra H^0(B,R^1f_*(\O_V))\lra H^2(B,\O_B),
\end{equation} 
and the term on the right vanishes. 
The following result reveals that the prime $p=2$ plays a special role in the general theory of linear systems. It
was already used in the previous section, and constitutes our second main result:

\begin{theorem}
\mylabel{pencils characteristic two}
In the above setting, suppose   that  $H^1(X,\O_X)=0$,
and that for each non-zero linear combination $t=\lambda_0s_0+\lambda_1s_1$, 
the resulting effective Cartier divisor $D_t\subset X$ is reduced and geometrically connected.
The the following holds:
\begin{enumerate}
\item The Stein factorization $B$ is a regular genus-zero curve.
\item The map $g:B\ra\PP^1$ is a universal homeomorphism, with $\deg(g)\mid 2$.
\item For each non-zero $s\in H^0(X,\shL)$ the resulting $D\subset X$ has  $h^0(\O_D)=\deg(g)$.
\item If $\deg(g)=2$ then the ground field $F$ is imperfect of characteristic $p=2$, and $B$ is a twisted line or twisted ribbon.
\end{enumerate}
\end{theorem}

\proof
Assertion (i)  immediately follows from the exact sequence \eqref{five term exact sequence}.
Given a rational point $t\in\PP^1$, we set $T=g^{-1}(t)$ and consider the schematic fiber $V_t=h^{-1}(t)=f^{-1}(T)$.
The projection $f:V_t=f^{-1}(T)\ra T$ is surjective and flat. 
Note that    $V_t$ is identified with the  effective Cartier divisor $D_t\subset X$ 
defined by the global section $t=\lambda_0s_0+\lambda_1s_1$. 
By assumption, $V_t$ is reduced and geometrically connected,
so the same holds for $T$. In particular, $T=\{b\}$ is a singleton, with coordinate ring $\kappa(b)$.
So Lemma \ref{radical degree two} below gives (ii) and (iv).  
 
It remains to verify (iii). 
Since the morphism $f:V\ra B$ is flat, the function $b\mapsto \dim_{\kappa(b)} H^0(V_b,\O_{X_b})$ is upper semicontinuous
(\cite{Hartshorne 1977}, Chapter III, Theorem 12.8). So the set $U\subset B$ where it takes the generic value
is open. If the residue field of $t=g(b)$ is separable,
the inclusion  $\Spec\kappa(b)\subset g^{-1}(t)$
is an equality, and hence
\begin{equation}
\label{dimension conversion}
\dim_{\kappa(b)} H^0(V_b,\O_{V_b}) = \deg(g)^{-1}\cdot \dim_{\kappa(t)} H^0(V_t,\O_{V_t}).
\end{equation}
Applying Lemma \ref{independent degrees} with the base-change $X\otimes\kappa(t)$, we see that
the  right-hand side does not depend on the separable point $t\in\PP^1$, and infer 
that $b\in U$ whenever the image  $f(b)\in\PP^1$ is separable.

By Grauert's Criterion (\cite{Hartshorne 1977}, Chapter III, Corollary 12.9), the formation of $f_*(\O_X)$ commutes with base-change over $U$. 
Thus $H^0(X_b,\O_{X_b})$ are one-dimensional vector spaces over $\kappa(b)$ for all $b\in U$.
If furthermore $t=f(b)$ is separable we get $h^0(\O_{X_t})=\deg(g)$ from \eqref{dimension conversion}. 
Now Lemma \ref{independent degrees} yields (iii).
\qed
 
\medskip
The above arguments rest  on the following key observation: 

\begin{lemma}
\mylabel{radical degree two}
Let $B$ be a regular genus-zero curve, and   $g:B\ra\PP^1$ be a surjective morphism  
of   degree $d\geq 2$. Suppose    that   
for each rational point $t\in\PP^1$, the fiber $g^{-1}(t)$ is reduced and geometrically connected.
Then   the field $F$ is imperfect of characteristic $p=2$,  the curve $B$ is a twisted line or twisted ribbon,
and $g:B\ra\PP^1$ is a universal homeomorphism of degree $d=2$.
\end{lemma}

\proof
Suppose the field $F$ is perfect.  Choose a geometric point $\Spec(\Omega)\ra \PP^1$ over some
rational point $t\in\PP^1$. The resulting geometric fiber is both reduced and connected, hence
isomorphic to $\Spec(\Omega)$, which results in the contradiction $d=1$.
Thus $F$ is imperfect, and we are in characteristic $p>0$.

Let $K=\O_{\PP^1,\eta}$ and $L=\O_{B,\eta}$ be the function fields of our curves. We start with an observation on   intermediate
fields $K\subsetneqq L'\subset L$. 
In light of \cite{EGA II}, Proposition 7.4.18 and Corollary 7.4.13, each such $L'$ 
defines  a regular curve $B'$, together with a factorization $g=g'\circ h$
into $h:B\ra B'$ and $g':B'\ra\PP^1$. Note that all three morphisms are   surjective, finite and flat.
We claim that $B'(F)=\varnothing$. Indeed, if there is a rational point $b'\in B'$, the image $t=g'(b')$ is rational,
and the fiber $Z=h^{-1}(b)$ is  a finite subscheme with $h^0(\O_Z)=[L':K]\geq 2$. Since $b\in Z$, 
the scheme $Z$ is non-reduced or disconnected, in contradiction to our assumption on the fibers of $g:B\ra\PP^1$.

Applying this observation with $L'=L$, we see that $B$ contains no rational point. 
Let  $L'$ be the relative separable closure of $K\subset L$.
Then the generic fiber of $g':B'\ra \PP^1$ is \'etale, so there is a non-empty open set $U\subset \PP^1$ over which $g'$ is \'etale.
Since $F$ is infinite, there must be rational points
$t\in U$. If  $[L':K]=\deg(B'/\PP^1)$
is greater than one, the fiber $Z=g^{-1}(t)$ is   geometrically disconnected,   contradiction.
It follows that
the field extension $K\subset L$ is purely inseparable (\cite{A 4-7}, Chapter V, \S7, No.\ 8, Proposition 13),
hence  $g:B\ra\PP^1$ is a universal homeomorphism, of degree $d=p^r$ for some $r\geq 1$.
If $p\neq 2$   the invertible sheaf $g^*\O_{\PP^1}(1)$ has odd degree, in contradiction to Proposition \ref{classification  genus-zero curves}.
This establishes $p=2$. From $B(F)=\varnothing$ we  also see that  $B$ is either a twisted line or a twisted ribbon.

Since $K\subset L$ is purely inseparable, we find a chain of   intermediate fields $K=L_0\subset \ldots\subset L_r=L$
with $[L_{i+1}:L_i]=p$. Consider the corresponding  regular curves $B_i$ and finite surjective morphisms
$$
B=B_r\lra B_{r-1}\lra \ldots\lra B_0=\PP^1.
$$
It remains to show $r=1$, and for this it suffices to verify that $B'=B_{r-1}$ is a genus-zero curve containing a rational point.
Let $f:B\ra B'$ be the given morphism. The short exact sequence
\begin{equation}
\label{tschirnhausen sequence}
0\lra \O_{B'}\lra f_*(\O_B)\lra \shL^{\otimes-1}\lra 0,
\end{equation}
defines an invertible sheaf $\shL$ on $B'$. The trace map $f_*(\O_B)\ra\O_{B'}$ vanishes on the subsheaf
$\O_{B'}$ because $d=p>0$, and the induced map $\shL^{\otimes-1}\ra\O_{B'}$ is non-zero since $B$ is reduced.
Hence $\deg(\shL)\geq 0$. On  the other hand, the long exact sequence for \eqref{tschirnhausen sequence}
reveals
$$
h^0(\O_{B'})=1\quadand h^0(\shL^{\otimes-1})= h^1(\O_{B'})\quadand h^1(\shL^{\otimes-1})=0.
$$
This gives $\deg(\shL)= \chi(\O_{B'})- \chi(\shL^{\otimes-1})=1-2h^1(\O_{B'})$. 
Consequently $B'$ is a genus-zero curve and  $\shL$ has degree one. It follows that $B'$ contains rational points.
\qed

\section{Bend and break}
\mylabel{Bend and break}

We keep the notation as in Section \ref{Surfaces without}, such that $X$ is a minimal regular surface  
with invariants $h^0(\O_X)=1$ and $h^1(\O_X)=h^0(\omega_X^{\otimes2})=0$, over a ground field $F$. 
The third  main results of this paper is the following:

\begin{theorem}
\mylabel{not integral}
Suppose there is an integral curve $C\subset X$ with $h^1(\O_C)=0$ and $C^2=4$. Then the curve is linearly equivalent to some $C'\subset X$ that is not integral.
\end{theorem}

In other words, ``bending'' the curve $C$ within its linear system, it ``breaks'' into some $C'=C'_1+C'_2$.
The proof requires  extensive preparation, and will be completed in Section \ref{Proofs}.
Special cases will already treated at the respective ends of this and the following two sections.

Our approach depends on certain maps $f:X\ra\PP^3$ that we introduce now.
First note that $h^0(\O_C)=1$, according to Proposition \ref{field of constants}.  
By Proposition \ref{globally generated}, the  invertible sheaf $\shL=\O_X(C)$ is globally generated, the restriction map $H^0(X,\shL)\ra H^0(C,\shL|C)$
is surjective, 
with $h^0(\shL)= 6$ and   $h^0(\shL|C)=5$.

We now choose a four-dimensional linear system as follows: 
Select two   sections on $\shL|C$ without common zeros, and extend them to global sections $s_1,s_2\in H^0(X,\shL)$.
Choose one more global section $s_3$
so that the restrictions $s_1|C,s_2|C,s_3|C$ are linearly independent, and let $s_0$ be a global section
whose zero-scheme is $C$. This defines a morphism 
$$
f:X\lra \PP^3\quad\text{with}\quad f^*\O_{\PP^3}(1)=\shL,
$$
having  $f^{-1}(t_i)\otimes 1=s_i$ when writing $\PP^3=\Proj F[t_0,\ldots,t_3]$. 
The schematic image $V\subset \PP^3$ is an integral surface. It is not a linear subscheme,
because the canonical  map $H^0(\PP^3,\O_{\PP^3}(1))\ra H^0(X,\shL)$ is injective.

The induced morphism $f:X\ra V$ is an alteration.
Let $Y=\Spec f_*(\O_X)$ be the Stein factorization, and write 
$$
h:X\lra Y\quadand g:Y\lra V
$$
for the resulting morphisms. 

\begin{proposition}
\mylabel{flatness of maps}
The finite morphism $g:Y\ra V$ is flat over the regular locus $\Reg(V)$,
and the  birational morphism $h:X\ra Y$ is the minimal resolution of singularities
for the normal surface $Y$.
\end{proposition}

\proof
The scheme $Y$ is Cohen--Macaulay, and over
the complement of $\Sing(V)$ the finite surjective  morphisms $g:Y\ra V$ must be flat
(\cite{Serre 1965},  page IV-37, Proposition 22).  The second statement holds because the regular surface $X$ is minimal. 
\qed

\medskip
To simplify notation, the pullbacks of the invertible sheaf $\O_{\PP^3}(1)$
are also written as  $\O_V(1)$ and $\O_Y(1)$ and $\O_X(1)=\shL$.

\begin{proposition}
\mylabel{image surface}
The scheme  $V$ is Gorenstein with $\omega_V=\O_V(d-4)$, where $d=\deg(V)$, and the numerical invariants
are $h^0(\O_V)=1$ and $h^1(\O_V)=0$. Moreover, the canonical maps  
\begin{equation}
\label{maps linear systems}
H^0(\O_{\PP^3}(1))\ra H^0(\O_V(1))\quadand H^0(\O_Y(1))\ra H^0(\O_X(1))=H^0(X,\shL)
\end{equation} 
are    bijective, and  $H^0(\O_V(1))\ra H^0(\O_Y(1))$ is injective.
\end{proposition}

\proof
Since $V\subset\PP^3$ is an effective Cartier divisor, it must be Gorenstein, and the Adjunction Formula gives
$\omega_V=\O_V(d-4)$. From the long  exact sequence stemming from  $0\ra \O_{\PP^3}(-d)\ra \O_{\PP^3}\ra\O_V\ra 0$,
we immediately get the numerical invariants.

From the long exact sequence for 
$0\ra \O_{\PP^3}(1-d)\ra \O_{\PP^3}(1)\ra \O_V(1)\ra 0$ and the vanishing of $H^i(\PP^3,\O_{\PP^3}(1-d))$ for $i\leq 1$
one sees that the restriction  map $H^0(\PP^3,\O_{\PP^3}(1))\ra H^0(V,\O_V(1))$ is bijective.
We have  $h_*(\O_X)=\O_Y$, so the Projection Formula ensures that $H^0(\O_Y(1))\ra H^0(\O_X(1))$ is bijective.
By construction, the composite map $H^0(\O_{\PP^3}(1))\ra H^0(X,\shL)$ is injective, hence the same holds for
$H^0(\O_V(1))\ra H^0(\O_Y(1))$.
\qed

\medskip
From  \eqref{degree equation} we see that there are 
two cases: Either $V\subset \PP^3$ is a quadric surface and $g:Y\ra X$ is a double covering,
or $V\subset \PP^3$ is a quartic surface  and $g:Y\ra V$ is birational.
The former case, together with other special situations,  can be treated quickly:

\begin{proposition}
\mylabel{thm in easy cases}
Theorem \ref{not integral} holds under any of the  following conditions:
\begin{enumerate}
\item There is a closed point $x\in X$ of degree at least five mapping to a   point $v\in V$ of degree at most two.
\item The singular locus of $Y$ contains a   point  whose image on $V$ has degree at most three.
\item The normal surface $Y$ contains a closed point  of degree three.
\item The surface $V\subset\PP^3$ has degree two.
\end{enumerate}
\end{proposition}

\proof
(i) Suppose some closed  point $x\in X$ has $[\kappa(x):F]\geq 5$ and maps  to a  point $v\in V$ with $[\kappa(v):F]\leq 2$.
Then there are  two planes $H_1\neq H_2$ with $v\in H_1\cap H_2$. Set   $C_i=f^{-1}(D_i)$.
We are done one of the $C_i$  non-integral, so we assume that both are integral. Then $C_1\cap C_2$ is zero-dimensional,
with $h^0(\O_{C_1\cap C_2})=(C_1\cdot C_2)=C^2=4$. By construction $x\in C_1\cap C_2$,
giving $4\geq h^0(\O_{C_1\cap C_2,x})\geq 5$, contradiction.
 
(ii) Suppose some $y\in \Sing(Y)$ has an image $v=g(y)$ with $[\kappa(v):F]\leq 3$.
Using $H^0(\PP^3,\O_{\PP^3}(1))=4$, we find some plane $H\subset \PP^3$ containing $v$.
Then $C'=f^{-1}(H)$ is a curve linearly equivalent to $C\subset X$. It contains the strict transform of the
curve $H\cap V$, together with some exceptional divisor, hence is reducible.

(iii) Now we have a closed point $y\in Y$ of degree three. In light of the previous cases,  we merely have to treat the case
that the local ring $\O_{Y,y}$ is regular.  Again we find a plane $H\subset\PP^3$ containing $v$.
Then $C'=f^{-1}(H)$ is a curve linearly equivalent to $C$, containing a point of degree three.
It is thus  a copy of the projective line, and we find two rational points $x\neq x'$ on $X$. Let $Z\subset X$
be the closed subscheme of degree four with coordinate ring $\O_{X,x}/\maxid_x^2\times\kappa(x')$, and recall that  $\shL=\O_X(C)$ has $h^0(\shL)=C^2+2=6$. 
So some curve $C'$ that is linearly equivalent to $C$ contains $Z$. By  Proposition \ref{classification genus-zero curves}, 
this curve is non-integral.

(iv) Now  the finite morphism  $g:Y\ra V$ has degree two.
Seeking a contradiction, we assume that every curve $C'\subset X$ linearly equivalent to $C$ is integral.
The   alteration $f:X\ra V$ factors over the normalization $\tilde{V}$, hence
$h^0(\O_{\tilde{V}})\leq h^0(\O_X)=1$. Proposition \ref{connectedness and normalization} ensures that $V$ is normal,
and that the singular locus contains at most one point. 
For each plane $H\subset \PP^3$, the   intersection $D=V\cap H$ and its preimage $C'=f^{-1}(D)$ are integral genus-zero curves.
The latter is linearly equivalent to $C$,  and the morphism $f:C'\ra D$ has  degree two.

Now choose $H$ disjoint from $\Sing(V)$. Then $f:C'\ra D$ is flat, and we see from Proposition \ref{covering implies line}
that $D$ is isomorphic to $\PP^1$. It follows that the regular locus of $V$ contains rational points. Fix such a rational point $v\in V$
and set $Z=\Spec(\O_{V,v}/\maxid^2)$. Using $h^0(\O_Z)=3$ and $h^0(\O_{\PP^3}(1))=4$ we can find another plane $H\subset \PP^3$,
now containing $Z$. This gives new integral genus-zero curves $D=V\cap H$ and $C'=f^{-1}(D)$, related by a morphism
$f:C'\ra D$ of degree two. By construction, the local ring $\O_{D,v}$ is singular.  According to Lemma \ref{connectedness and normalization},
our plane $H$ must be disjoint from $\Sing(V)$, and it follows that $f:C'\ra D$ is flat. As above we get 
$D=\PP^1$. But this contradicts that $D$ is singular.
\qed

\section{Non-normal quartic surfaces}
\mylabel{Non-normal quartic}

We keep all assumptions of the preceding section, and continue to work towards the proof for Theorem \ref{not integral}.
Recall that $X$ is a  minimal regular surface  over the ground field $F$
with invariants $h^0(\O_X)=1$ and $h^1(\O_X)=h^0(\omega_X^{\otimes2})=0$, containing an  integral curve $C\subset X$
with $h^1(\O_C)=0$ and $C^2=4$. 
\emph{We now additionally assume that our linear system maps $X$   
to a quartic surface $V\subset\PP^3$.} 
Then    $f:X\ra V$ is birational,
$g:Y\ra V$ is the normalization, and $h:X\ra Y$ is the minimal resolution of singularities.
 
\begin{proposition}
\mylabel{non-normal}
The quartic surface $V\subset\PP^3$ is non-normal.
\end{proposition}

\proof
Suppose $V$ is normal, such that $Y=V$. Form the transcendental field extension $F'=F(t)$.
One easily checks that $X\otimes F'$ remains regular and minimal, and $C\otimes F'$ stays integral, and $V$ is normal if and only if this holds
for $V\otimes F'$. So without loss of generality we may assume that 
$F$ is infinite.
We then  find a plane $H\subset\PP^3$ that avoids the finite set $\Sing(V)$.
Then the projection of $C'=f^{-1}(D)$ to $D=H\cap V$ is an isomorphism.
By the Genus Formula, the  quartic plane curve  $D\subset H$ has invariants $h^0(\O_D)=1$ and $h^1(\O_D)=(4-1)(4-2)/2=3$.
On the other hand, $C'\subset X$ is linearly equivalent to $C$, hence has $h^1(\O_{C'})=0$, contradiction.
\qed

\medskip
The normalization map $g:Y\ra V$ comes with  a \emph{branch curve} $B\subset V$, defined as 
the schematic support of $g_*(\O_Y)/\O_V$. Its schematic preimage $R=g^{-1}(B)$ is called the
\emph{ramification curve}. 
Both are indeed equidimensional of dimension one, and  without embedded components, because
$Y$ and $V$ are Cohen--Macaulay.
Relative Duality  for the normalization map gives the formula $g_*(\omega_{Y/V})= \uHom(g_*(\O_Y),\O_V)$.
The latter coincides with $g_*(\O_Y(-R))$,  and we obtain
\begin{equation}
\label{adjunction formula}
\omega_Y=\omega_{Y/V}\otimes f^*(\omega_V)=\omega_Y=\O_Y(-R).
\end{equation} 
For more details we refer to \cite{Fanelli; Schroeer 2020a}, Appendix A.
The situation is summarized in the   commutative diagram
\begin{equation}
\label{conductor square quartic}
\begin{tikzcd} 
		& X\ar[d,"h"']\ar[dd, bend left=40, "f"]\\
R\ar[r]\ar[d]	& Y\ar[d,"g"']\\
B\ar[r]		& V\ar[r]		& \PP^3,
\end{tikzcd}
\end{equation}
all horizontal arrows are closed embeddings. Also note that the annihilator ideal for $g_*(\O_Y)/\O_V$ is called
the \emph{conductor ideal} for $g:Y\ra V$, and accordingly the  lower left part of the diagram is called the \emph{conductor square}.

In our situation,  the branch curve $B$ can also be seen as a  \emph{space curve} in $\PP^3$, which will be crucial throughout.
Recall that each space curve $Z\subset\PP^3$ has some \emph{degree}, defined as 
$$
\deg(Z) = \chi(\O_Z(1)) - \chi(\O_Z) = \deg(\O_Z(1)) = (\O_{\PP^3}(1)\cdot C)   \geq 1.
$$
Space curves  $Z$ of degree one are called  \emph{lines}. They are  integral, with   $\O_Z(1)$  globally generated of degree one. This sheaf is then 
generated by two global sections, and the resulting morphism $Z\ra\PP^1$ of degree one is an isomorphism.
One finds that the lines   are precisely the intersections of two  different planes.
Note that space curves $Z$ contained in a plane  are rather special: By the Genus Formula they have
invariants $h^0(\O_Z)=1$ and $h^1(\O_Z)=(d-1)(d-2)/2$, where $d=\deg(Z)$.

\begin{proposition}
\mylabel{degree branch curve}
The branch curve $B\subset \PP^3$ has degree three. Moreover, $B$ is integral provided that it contains no line.
\end{proposition}

\proof
Let $\eta\in B$ be some generic point, and set $\Lambda=\O_{B,\eta}$ and $\Lambda'=f_*(\O_R)_\eta$.
According to \cite{Fanelli; Schroeer 2020a}, Proposition  A.2 the lengths of these $\Lambda$-modules are related by  $\length (\Lambda')=2\length(\Lambda)$.
In light of \cite{Kleiman 1966},  Proposition  6 on page 299, this ensures $\deg(\shN_R) = 2\deg(\shN)$ for
every invertible sheaf $\shN$ on the branch curve.

Let $E_1,\ldots,E_r$ be the irreducible components of the exceptional divisor $\Exc(X/Y)$,
endowed with reduced scheme structure, and $R'\subset X$ be the strict transform of the ramification curve
$R\subset Y$. Write
$$
K_{X/Y} = -\sum\lambda_i E_i\quadand f^*(R) = R'+\sum\mu_iE_i
$$
for certain  $\lambda_i,\mu_i\in\QQ$. Then the numerical class of $\omega_X$ is given by the
$\QQ$-divisor $K_{X/Y} - f^*(R)$. With $\shL=\O_V(1)$ we compute
$$
(\omega_X\cdot C) = (K_{X/Y}- f^*(R))\cdot f^*(\shL) = - R\cdot g^*(\shL) = -2\deg(\shL|C) = -2 \deg(B),
$$
using the Projection Formulas. (In the above arguments  we have use Mumford's rational pullbacks and rational intersection numbers;
for details see  the discussion in \cite{Mumford 1961} and \cite{Schroeer 2019}).
On the other hand, the Adjunction Formula gives 
$$
(\omega_X\cdot C) = \deg(\omega_C) - C^2= -2  - 4 =-6.
$$
Combining the above equations yields $\deg(B) = 3$. 

Suppose that   $B$ contains no line.
Let $\eta_1,\ldots,\eta_r\in B$ be the generic points and $B_1,\ldots,B_r\subset B$ 
be the corresponding  irreducible components, viewed as the schematic closure of the canonical morphism $\Spec(\O_{B,\eta_i})\ra B$.
Then 
$$
3=\deg(B)=\sum_{i=1}^r\deg(B_i) =  \sum_{i=1}^r \deg(B_{i,\red})\cdot \length(\O_{B_i,\eta_i})\geq 2lr,
$$
where $l\geq 1$ is the smallest among the $\length(\O_{B_i,\eta_i})$.
We see $r=1$ and $l=1$, hence  $B$ is integral.
\qed

\medskip
The presence of lines is of little interest:

\begin{proposition}
Theorem \ref{not integral} holds if the branch curve $B$ contains a line.
\end{proposition}

\proof
Choose a line $L\subset B$, and two planes $H_1\neq H_2$ inside $\PP^3$ containing the line.
If $H_i\cap V$ is reducible, the same holds for the preimage $C_i=f^{-1}(H_i\cap C)$.
We thus may assume that $H_1\cap V$ and $H_2\cap V$ have the same support,
namely $L$. Then the  preimages $f^{-1}(H_i\cap V)$ are   effective Cartier divisors
on some normal scheme that are linearly equivalent and have the same support. This is only possible if
$f^{-1}(H_1\cap V)=f^{-1}(H_2\cap V)$. But this contradicts the injectivity
of the map $H^0(\PP^3,\O_{\PP^3}(1))\ra H^0(X,\shL)$, where $\shL=\O_X(C)$.
\qed

\section{Twisted  cubics and exotic cubics}
\mylabel{Twisted cubics}
 
We keep the assumptions from the preceding section, and continue work towards the proof for Theorem \ref{not integral}.
Recall that our regular surface $X$   over the ground field $F$ contains an  integral curve $C\subset X$
with $h^1(\O_C)=0$ and $C^2=4$, and maps to non-normal quartic surface $V\subset\PP^3$,
with normalization $g:Y\ra V$, and resulting branch curve $B\subset V$ and ramification curve $R\subset Y$. 
\emph{We now additionally assume that  $B$ contains no line, and that $\Sing(Y)$ contains no point whose image on $V$ is rational.}

These two innocuous assumptions are of profound consequence. Among other things, they allow
to compute the cohomological invariants of our schemes:
 
\begin{proposition}
\mylabel{invariants surfaces curves}
The  normal surface $Y$ has at most rational singularities,
and the numerical invariants are $h^0(\O_Y)=1$ and $h^1(\O_Y)=h^2(\O_Y)=0$. Moreover, the ramification curve
and the branch curve have invariants
$$
h^0(\O_R)=h^1(\O_R)=1\quadand h^0(\O_B)=1,\, h^1(\O_B)=0,
$$
and the scheme $B$ is integral.
\end{proposition}

\proof
First note that $B$ is integral by Proposition \ref{degree branch curve}.
The Leray--Serre spectral sequence for  $h:X\ra Y$   gives $h^0(\O_Y)=1$ and $h^1(\O_Y)=0$,
together with an identification $H^2(Y,\O_Y)=H^0(Y,R^1h_*(\O_X))$. 
The conductor square in \eqref{conductor square quartic} yields a short exact sequence $0\ra\O_V\ra\O_Y\oplus \O_B\ra\O_R\ra 0$,
and we proceed by examining its long exact sequence.  
It starts with   
$$
0\lra H^0(\O_V)\lra H^0(\O_Y)\oplus H^0(\O_B)\lra H^0(\O_R)\lra H^1(\O_V).
$$
The term on the right vanishes, and we see $h^0(\O_B)=h^0(\O_R)$.
Using $h^1(\O_V)=h^1(\O_Y)=0$, we get another    exact sequence
$$
0\lra H^1(\O_B)\lra H^1(\O_R)\lra H^2(\O_V)\lra H^2(\O_Y)\lra 0.
$$
Suppose   that the map on the right is non-zero. It is necessarily bijective, because $h^0(\O_V)=1$,  and thus $h^2(\O_Y)=1$.
In turn, the sheaf $R^1h_*(\O_X)$ must be supported by a rational point $y\in Y$. The image $v\in V$ of this singularity
is again a rational point, in contradiction to our standing assumption.
Consequently $h^2(\O_Y)=0$ and $R^1h_*(\O_X)=0$, such that each singular local ring $\O_{Y,y}$
is a rational singularity, having some non-trivial field extension $F\subset\kappa(y)$. Moreover, we see that $h^1(\O_R)=h^1(\O_B)+1$.

The short exact sequence $0\ra\O_Y(-R)\ra \O_Y\ra\O_R\ra 0$ gives a long exact sequence
$$
H^0(Y,\O_Y)\lra H^0(R,\O_R)\lra H^1(Y,\O_Y(-R)).
$$
Using $\omega_Y=\O_Y(-R)$ and Serre Duality, we see that the term on the right is 
is dual to $H^1(Y,\O_Y)$, which vanishes. Thus $h^0(\O_R)=h^0(\O_Y)=1$. Finally, we have an exact sequence
$$
H^1(Y,\O_Y)\lra H^1(R,\O_R)\lra H^2(Y,\O_Y(-R))\lra H^2(Y,\O_Y).
$$
The outer terms vanish, so $h^1(\O_R)=h^2(\O_Y(-R))=1$, again  by Serre Duality.
\qed

\medskip
A useful geometric consequence:

\begin{proposition}
\mylabel{factorization over blowing-up}
The resolution of singularities $f:X\ra V$ factors over the blowing-up $\Bl_B(V)\ra V$ with respect to the space curve $B\subset V$.  
\end{proposition}

\proof
According to Hartshorne \cite{Hartshorne 1977}, Chapter II, Proposition 7.14
the task is  verify that $f^{-1}(Z)\subset X$ is an effective Cartier divisor.
As  $g^{-1}(Z)=R$,   it suffices to show that $h^{-1}(R)\subset X$ is an effective Cartier divisor.
Since $Y$ has only rational singularities, this indeed holds, according to \cite{Schroeer 2007b}, Proposition 10.5.
\qed

\medskip
Recall that a space curve $Z\subset\PP^3$ of degree three that is isomorphic to the projective line is called a \emph{twisted cubic}.
The Adjunction Formula ensures that it is not contained in a plane.
Furthermore   $\O_Z(1)\simeq\O_{\PP^1}(3)$,  and  
the map $H^0(\PP^3,\O_{\PP^3}(1))\ra H^0(Z,\O_Z(1))$ is injective, hence bijective.
In turn, each twisted cubic can be seen as the third Veronese embedding of $\PP^1$.
Moreover, the automorphism group $\PGL_4(F)$ acts transitively on the set of twisted cubics.
Note that that  the classical designation ``twisted''  refers to the tensor power $\O_{\PP^1}(3)$, and has nothing to do with
modern usage in connection with  twisted forms.

We also have  to deal with a more outlandish type of integral space curve  of degree three,
which  are of arithmetic nature and depend on some  cubic field extension $F\subset E$.
Use such an  extension to form the cocartesian square
$$
\begin{CD}
\Spec(E)	@>(1:0)>> 	\PP^1_E\\
@VVV		@VVV\\
\Spec(F)	@>>>	Z.
\end{CD}
$$
Then $Z$ is a genus-zero curve that is non-Gorenstein.
One easily checks that $\Pic(Z)\ra\Pic(\PP^1_E)$ is bijective,
and that the invertible sheaf $\shL$ with $\shL|\PP^1_E=\O_{\PP^1}(1)$ is very ample, having $\deg(\shL)=3$ and   $h^0(\shL)=4$.
This yields a closed embedding $Z\subset \PP^3$ of degree three, where we suppress the dependence on the cubic field
extension   from notation. Note that this are cones over any embedding $\Spec(E)\subset\PP^2$.
For lack of better designation, we  call such space  curves $Z\subset\PP^3$  
\emph{exotic cubics}. Since the Hilbert scheme of space curves of degree three and genus zero is   irreducible
(\cite{Martin-Deschamps; Perrin 1996}, Theorem 4.1), one may view the exotic cubics as arithmetic degenerations of twisted cubics.

\begin{proposition}
\mylabel{integral branch curve}
Let $Z\subset \PP^3$ be an integral space curve of degree three, with invariants $h^0(\O_Z)=1$ and $h^1(\O_Z)=0$.
Then $Z$ is   a twisted cubic   or an exotic cubic.
\end{proposition}

\proof
First note that $\O_Z(1)$ is an invertible sheaf of odd degree. So if the genus-zero curve $Z$ is Gorenstein,
it must by isomorphic to the projective line, by Proposition \ref{classification genus-zero curves}, and hence is a twisted cubic.
Suppose now that $Z$ fails to be Gorenstein. The normalization $\tilde{Z}$ has $h^1(\O_{\tilde{Z}})=0$,   
the pullback $\O_{\tilde{Z}}(1)$ has degree three, and we infer that $\tilde{Z}=\PP^1_E$
for some field extension $F\subset E$ whose degree divides three.
Let 
$$
\begin{CD}
\tilde{Z}'	@>>>	\tilde{Z}\\
@VVV		@VVV\\
Z'	@>>>	Z.
\end{CD}
$$
be the conductor square, and $z_1,\ldots,z_r\in Z$ be the branch points. Set
$$
d=[E:F]\quadand l_i=[\O_{\tilde{Z}',z_i}:E]\quadand \lambda_i= dl_i - [\O_{Z',z_i}:F].
$$
From the short exact sequence $0\ra \O_Z\ra\O_{\tilde{Z}}\times\O_{Z'}\ra \O_{\tilde{Z}'}\ra 0$ we get
$$
1-(d +  \sum_{i=1}^r(dl_i - \lambda_i) ) +  d\sum_{i=1}^r l_i = 0.
$$
This  arrangement of terms produces   a \emph{partition} $d-1=\sum_{i=1}^r \lambda_i$. From $r\geq 1$ and $\lambda_i\geq 1$ we 
infer $d=3$, and the only possibilities are
$$
r=1,\, \lambda_1=2\quadand r=2,\,\lambda_1=\lambda_2=1.
$$
Suppose $r=2$. Then   $\O_{Z',z_i}\subset \O_{\tilde{Z}',z_i}$ has $F$-codimension $\lambda_i=1$,
so the subfield  $E_i=E\cap \O_{Z',z_i}$ inside $E$ has $F$-codimension at most one. 
From $H^0(Z,\O_Z)=F$ and the exact sequence
$$
0\lra H^0(Z,\O_Z)\lra H^0(\tilde{Z},\O_{\tilde{Z}})\times H^0(Z',\O_{Z'})\lra H^0(\tilde{Z}',\O_{\tilde{Z}'})
$$
we infer that at least one inclusion $E_i\subset E$ is strict. However, $[E:F]=3$  precludes such subfields.

Summing up, there is a unique branch point $z\in Z$, and   $\O_{Z',z}\subset\O_{\tilde{Z}',z}$ has $F$-codimension two.
From $h^1(\O_Z)=0$ one immediately sees that the preimage of $z\in Z$ is a single point $\tilde{z}\in\tilde{Z}$,
and that the residue field extension $\kappa(z)\subset\kappa(\tilde{z})$ is non-trivial.
Its $F$-codimensions is at most two. This codimension equals the $F$-dimension of $\kappa(z)$-vector space $\kappa(\tilde{z})/\kappa(z)$,
hence $d=[\kappa(z):F]$ divides two.
Consider the commutative diagram
$$
\begin{CD}
0	@>>>	F	@>>>	E\\
@.	@VVV		@VVV\\
0	@>>>	\kappa(z)	@>>>	\kappa(\tilde{z})	@>>> 	\kappa(\tilde{z})/\kappa(z)	@>>> 0.
\end{CD}
$$
The case $d=2$ yields $[\kappa(\tilde{z}):F]=4$, which is also a multiple of $[E:F]=3$, contradiction.
Thus $\kappa(z)=F$, and it follows $\kappa(\tilde{z})=E$.
The canonical surjection $\O_{\tilde{Z}',z}/\O_{Z',z}\ra\kappa(\tilde{z})/\kappa(z)=E/F$ is bijective,
because both quotients are two-dimensional $F$-vector spaces.
By the very definition, the conductor ideal $\idealc\subset\O_{Z,z}$ is the annihilator ideal
of $\O_{\tilde{Z},z}/\O_{Z,z}= \O_{\tilde{Z}',z}/\O_{Z',z}$, which here coincides with $\kappa(\tilde{z})/\kappa(z)$.
Consequently $\O_{\tilde{Z}',z}=\kappa(\tilde{z})=E$ and $\O_{Z',z}=\kappa(z)=F$. In other words, $Z$ is an exotic cubic.
\qed

\begin{proposition}
\mylabel{thm for exotic cubics}
Theorem \ref{not integral} holds provided the branch curve   $B$ is an exotic cubic.
\end{proposition}

\proof
Write $v_0\in B$ for the unique singularity of the branch curve for the normalization $g:Y\ra V$, which is a rational point.
The strategy  is to verify that the ramification curve $R=g^{-1}(B)$ contains
a closed point of degree three, or that all points in $g^{-1}(v_0)$ have degree at least five.
Then Proposition \ref{thm in easy cases} indeed gives our assertion.
To carry out this strategy  we proceed as follows:
The  normalization   takes the form $\tilde{B}=\PP^1_E$ for some cubic field extension $F\subset E$.
If $R$ is not integral, then both maps  $\tilde{B}\ra B\leftarrow R$ admit  sections over the complement
of the singularity $v_0\in B$, and $R$ indeed contains closed points of degree three.

The main task is to treat the case that $R$ is integral.
For each $y\in R$ mapping to  the   singularity $v_0\in B$
the local ring $\O_{Y,y}$ is regular,  by our standing assumption on the map  $\Sing(Y)\ra V$, and hence  $\O_{R,y}$ is Gorenstein.
Moreover, the morphism $g:R\ra B$ is flat of degree two on the complement of $v_0$, and we infer that
the whole curve $R$ is Gorenstein. 
Consider the commutative square
\begin{equation}
\label{branch curve and ramfication curve}
\begin{CD}
\tilde{B}	@<<<	\tilde{R}\\
@VVV		@VVV\\
B	@<<<	R
\end{CD}
\end{equation} 
where the vertical arrows are the normalizations, and the horizontal morphisms have degree two. The number
$h^1(\O_{\tilde{R}})$ is bounded above by $h^1(\O_R)=1$
and is also a multiple of $h^0(\O_{\tilde{B}})=3$, which is only possible when $h^1(\O_{\tilde{R}})=0$.
Thus $\tilde{R}$ is a genus-zero curve   over some further field extension $\tilde{E}=H^0(\tilde{R},\O_{\tilde{R}})$,
with  $E\subset\tilde{E}$ 
of degree either one or two.
We next form   the conductor square
$$
\begin{CD}
\tilde{R}'	@>>>	\tilde{R}\\
@VVV		@VVV\\
R'	@>>>	R.
\end{CD}
$$
for the normalization $\tilde{R}\ra R$. Write $y_1,\ldots,y_r\in R$ for  the branch points   and set
$$
d=[\tilde{E}:F]\quadand l_i=[\O_{\tilde{R}',y_i}:\tilde{E}]\quadand \lambda_i=dl_i-[\O_{R',y_i}:F].
$$
The short exact sequence $0\ra \O_R\ra \O_{\tilde{R}}\times\O_{R'}\ra \O_{\tilde{R}'}\ra 0$
yields
$$
1-(d + \sum_{i=1}^r(   dl_i-\lambda_i) ) + d\sum_{i=1}^r   l_i     -1=0,
$$
which gives a partition $d=\lambda_1+\ldots+\lambda_r$. We have  $\lambda_i=dl_i/2$ because $R$ is Gorenstein
(\cite{Fanelli; Schroeer 2020b}, Proposition A.2),
and thus get   another  partition $2=l_1+\ldots+l_r$.
The only possibilities are 
$$
r=1,\, l_1=2\quadand r=2,\, l_1=l_2=1.
$$
Moreover  $d=3$ or $d=6$. We now have to go through all possible cases:

Suppose first that $r=2$. In other words, the normalization produces exactly two branch points $y_1,y_2\in R$,
each having  $\O_{\tilde{R}',y_i}=\tilde{E}$. Inside this,  $\O_{R',y_i}$ is a subfield of index two.
We see  $\O_{R',y_i}=\kappa(y_i)$ and  the integer $d$ is even, therefore $d=6$, hence the closed points  $y_i\in Y$
have degree three, as desired.
 
Suppose now $r=1$. In other words, the normalization comes with exactly one branch point $y_1\in R$,
having $[\O_{\tilde{R}',y_1}:\tilde{E}]=2$. 
So the  $\tilde{E}$-algebra  $\O_{\tilde{R}',y_1}$ is either a quadratic field extension,  or the field product
$\tilde{E}\times\tilde{E}$, or the ring of dual numbers $\tilde{E}[\epsilon]=\tilde{E}\oplus\epsilon\tilde{E}$.
In the first two cases, the subring $\O_{R',y_1}$ is a field, and thus coincides with $\kappa(y_1)$.
In the last two cases, the  normal  curve  $\tilde{R}$ contains an $\tilde{E}$-valued point, and is thus
isomorphic to the projective line $\PP^1_{\tilde{E}}$.
So for $d=3$, we find in all three cases some closed point $y\in R$ of degree three.

Assume now $d=6$, so the $F$-algebra $\O_{\tilde{R}',y_1}$ has degree twelve.
Suppose first that the subring $\O_{R',y_1}$ is a field. It thus  coincides with $\kappa(y_1)$,
necessarily with  $[\kappa(y_1):F]=6$,
and thus every closed point $y\in R$ has degree at least six, as desired.
We finally come to  the most interesting  case where  $\O_{\tilde{R}',y_1}=\tilde{E}[\epsilon]$ is the ring of dual numbers,
and that the subring $\O_{R',y_1}$ is not a field.
It takes the form $K+\epsilon U$, for a \emph{field of representatives}  $K$
and some non-zero $K$-vector subspace $\epsilon U\subset \epsilon \tilde{E}$, satisfying
$$
6=[ \O_{R',y_1}:F]= [K:F](1+\dim_K(U)) \quadand \dim_K(U)=\edim(\O_{R',y_1}).
$$
Here $\edim(A)=\dim_\kappa(\maxid/\maxid^2)$ denotes the \emph{embedding dimension} of a local noetherian ring $(A,\maxid,\kappa)$.
By construction,  $[K:F]$ is a strict divisor of $[\tilde{E}:F]= 6$.
If $[K:F]=3$ the closed point $y_1\in Y$ has degree three, as desired.
If $[K:F]=1$ we have $\edim(\O_{Y,y_1})\geq \edim(\O_{R',y_1})=5$.
Then $y_1$ is a rational point in $\Sing(Y)$, in contradiction to our standing assumption that no such point exists. 

Finally, we have to rule out the case  $[K:F]=2$ and $\dim_K(U)=2$. Then the canonical surjection $E\otimes_F K\ra \tilde{E}$
is bijective, and we have an identification $\tilde{R}=\tilde{B}\otimes_E\tilde{E}=\tilde{B}\otimes_FK$.
Since all closed points $y\in R$ except $y_1$ have degree at least six,
it suffices  to treat the case that
$g^{-1}(v_0)=\{y_1\}$, in light of Proposition \ref{thm in easy cases}.
 Write $\tilde{v}_0\in \tilde{B}$ and $\tilde{y}_1\in \tilde{R}$ for the preimages. 
From \eqref{branch curve and ramfication curve} we obtain a commutative diagram of cotangent spaces
$$
\begin{CD}
\maxid_{\tilde{v}_0}/\maxid^2_{\tilde{v}_0}	@>>>	\maxid_{\tilde{y}_1}/\maxid^2_{\tilde{y}_1}\\
@A\id AA			@AAA\\
\maxid_{v_0}/\maxid^2_{v_0}	@>>>	\maxid_{y_1}/\maxid^2_{y_1}.
\end{CD}
$$
These vector spaces and linear maps take the explicit form
$$
\begin{CD}
\epsilon E	@>>>	\epsilon (E\otimes_FK)\\
@A\id AA			@AAA\\
\epsilon E	@>>>	\epsilon U.
\end{CD}
$$
Choose an $F$-basis $\epsilon\alpha_i\in \epsilon E$, $1\leq i\leq 3$.
Their images $\epsilon(\alpha_i\otimes 1)$ in $\epsilon (E\otimes_FK)$
stays $K$-linearly independent. However, the images in the two-dimensional $K$-vector space $\epsilon U$
must be $K$-linearly dependent, contradiction.
\qed

\section{Blowing-ups centered at twisted cubics }
\mylabel{Blowing-ups}
 
We keep the assumptions from the preceding section, and continue work towards the proof for Theorem \ref{not integral}.
Recall that our regular surface $X$   over the ground field $F$ contains an  integral curve $C\subset X$
with $h^1(\O_C)=0$ and $C^2=4$, and maps to non-normal quartic surface $V\subset\PP^3$,
with normalization $g:Y\ra V$, and resulting branch curve $B\subset V$ and ramification curve $R\subset Y$. 
\emph{We now   assume that  
 $B$ is a twisted cubic, and that $\Sing(Y)$ contains no point whose image on $V$ is rational.}

Choosing an identification $B=\PP^1$
together with  $\O_B(1)=\O_{\PP^1}(3)$,  we see that the inclusion of the twisted cubic is given by
$(x_0:x_1)\mapsto (x_0^3:x_0^2x_1:x_0x_1^2:x_1^3)$, and its ideal $\ideala\subset k[t_0,\ldots,t_3]$
is generated by the three   homogeneous polynomials
\begin{equation}
\label{generators twisted cubic}
f_0=t_0t_4-t_1t_2,\quadand f_1=t_0t_3-t_1^2,\quadand f_2=t_1t_3-t_2^2.
\end{equation}
Let $P\in k[t_0,\ldots,t_3]$ be a homogeneous polynomial of degree four defining the quartic surface $V\subset\PP^3$.
The inclusion $B\subset V$ translates into  $P\in \ideala$. One actually can say much more:

\begin{proposition}
\mylabel{defining quartic polynomial}
There is a  homogeneous polynomial $\Phi\in F[X_0,X_1,X_2]$ of degree two such that 
$P=\Phi(f_0,f_1,f_2)$. Moreover, such $\Phi$ is irreducible.
\end{proposition}

\proof
Since $B$ is reduced and $V$ is singular along it, we must have $P\in\ideala^2$. 
Now write $P=\sum g_{ijd} f_if_j$ for some homogeneous   $g_{ijd}\in F[t_0,\ldots,t_3]$ of degree $d\geq 0$.
The terms where  $\deg(g_{ijd})+\deg(f_i)+\deg(f_j)$ differs from $\deg(P)=4$ cancel each other, and discarding them
we may assume that each non-zero summand  $d=0$. 
Summing up, we have  $P=\Phi(f_0,f_1,f_2)$   for some ternary quadratic form $\Phi$.
The latter must be irreducible, because $P$ is irreducible.
\qed

\medskip
To proceed, we consider  the blowing-up  $\Bl_B(\PP^3)\ra\PP^3$ with center the  twisted cubic $B$.
The strict transform of $V$ coincides with the blowing-up $\Bl_B(V)\ra V$.
According to Proposition \ref{factorization over blowing-up}, the morphism $f:X\ra V$ factors over $\Bl_B(V)$.
Let $\shI\subset\O_{\PP^3}$ be the sheaf of ideals for the twisted cubic.
The surjection
$$
\shH=\left(\bigoplus_{i=0}^2\O_{\PP^3}\right)\otimes\O_{\PP^3}(-2)=\bigoplus_{i=0}^2\O_{\PP^3}(-2)\stackrel{(f_0,f_1,f_2)}{\lra}\shI
$$
stemming form \eqref{generators twisted cubic} defines an inclusion $\Bl_B(\PP^3)\subset\PP(\shH)=\PP^3\times\PP^2$.
We thus obtain a commutative diagram
\begin{equation}
\label{ray blowing up}
\begin{tikzcd} 
X\ar[dd,"f"']\ar[rr] 	& 				& \Bl_B(V) \ar[d]\ar[ddll]\ar[ddrr]\\
 		&				& \Bl_B(\PP^3)\ar[dl,"\pr_1"']\ar[dr,"\pr_2"]	&	& \\
V\ar[r]		& \PP^3\ar[rr, dashed, "(f_0:f_1:f_2)"']	&				& \PP^2	& V_+(\Phi)\ar[l].\\
\end{tikzcd}
\end{equation} 
Indeed, the rational map $\PP^3\dashrightarrow \PP^2$ is defined outside $B$, and sends $V\smallsetminus B$ to the quadric curve $V_+(\Phi)$.
Using that $V\smallsetminus B$ is schematically dense in $\Bl_B(V)$,  
we immediately see that the induced projection $\pr_2:\Bl_B(V)\ra\PP^2$ factors over $V_+(\Phi)$, giving the   diagonal arrow to the right.

According to Ray's result (\cite{Ray 2020}, Corollary 3.6) the projection $\Bl_Z(\PP^3)\ra\PP^2$ is isomorphic to the projectivization $\PP(\shE)$,
with the locally free sheaf of rank two sitting in a short exact sequence
\begin{equation}
\label{sheaf giving projectivization}
0\lra \O_{\PP^2}^{\oplus 2}(-1) \lra \O_{\PP^2}^{\oplus 4}\lra \shE\lra 0.
\end{equation} 
Moreover, it follows from  loc.\ cit.\ Theorem 3.4 that 
\begin{equation}
\label{invertible sheaves on blowing-up}
\pr_2^*(\O_{\PP^2}(1)) \simeq \pr_1^*(\O_{\PP^3}(1))\otimes \O_{\Bl_B(\PP^3)}(-E),
\end{equation} 
where $E=\pr_1^{-1}(B)$ is the exceptional divisors. Note that Ray worked over the complex numbers, 
but his arguments literally hold true over arbitrary  ground fields.
In light of the commutative diagram \eqref{ray blowing up}, the morphism $X\ra\Bl_B(\PP^3)=\PP(\shE)$ factors over $\PP(\shE\mid V_+(\Phi))$.

\begin{theorem}
\mylabel{structure of blowing-up}
In the above situation, the following holds:
\begin{enumerate}
\item The integral quadric curve $V_+(\Phi)\subset\PP^2$ is regular.
\item The rank-two locally free sheaf  $\shE|V_+(\Phi)$   has degree four.
\item The morphism $X\ra  \PP(\shE|V_+(\Phi))$ is an isomorphism.
\item The curve $C\subset X$ is a section for the projection $X\ra V_+(\Phi)$.
\end{enumerate}
\end{theorem}

\proof
The projection $X\ra V_+(\Phi)$ is surjective, because the fibers for the bundle   $\pr_2:\Bl_B(\PP^3)\ra\PP^2$ are one-dimensional.
The projection thus factors  over the normalization of $V_+(\Phi)$. This  involves a constant field extension 
if $V_+(\Phi)$ is non-normal, in contradiction to  $H^0(X,\O_X)=F$.  From the exact sequence
\eqref{sheaf giving projectivization} we get $\det(\shE)=\O_{\PP^2}(2)$.
Since $V_+(\Phi)\subset\PP^2$ has degree two, the restriction $\shE \mid V_+(\Phi)$ has degree four. This establishes (i) and (ii).

The morphism $X\ra \PP(\shE|V_+(\Phi))$ is surjective, because  $\O_X(C)=f^*(\O_{\PP^3}(1))$ and $C^2>0$.
The scheme  $\PP(\shE\mid V_+(\Phi))$ is regular  because the base $V_+(\Phi)$ is regular, and (iii) follows from the minimality of our surface $X$.

It remains to establish (iv), which is the most interesting part. Our task is to verify that
the intersection number  $(\pr_2^*(\O_{\PP^2}(1))\cdot C)$ coincides
with $(\O_{\PP^2}(1)\cdot V_+(\Phi))=2$. 
To achieve this we compute  separately with the factors appearing in \eqref{invertible sheaves on blowing-up}:
The first factor contributes
$(\pr_1^*(\O_{\PP^3}(2))\cdot C) = 2\cdot C^2=8$.
The second factor is the negative of 
$$
(\O_{\Bl_B(\PP^3)}(E)\cdot C) = h^0(\O_{C\cap E})=h^0(\O_{C\cap X\cap E}) = \deg(\O_{\PP^3}(1)\mid R).
$$
But $\deg(\O_{\PP^3}(1)\mid R) =  \deg(R/B)\cdot \deg(B) = 2\cdot 3=6$, which gives the desired  $(\pr_2^*(\O_{\PP^2}(1))\cdot C)=8-6=2$.
\qed

\medskip
In particular, $X$ comes with a fibration such that the base and the generic fiber are genus-zero curves.
As discussed in Section \ref{Proofs}, this actually finalizes the proof for  
our generalization of Iskoviskih's result in Theorem \ref{trichotomy}.
To complete the proof for  the bend-and-break statement in 
Theorem \ref{not integral}, it remains to understand the geometry of ruled surfaces 
over a regular genus-zero curve, which boils down to an analysis
of locally free sheaves of rank two on on such curves.  We shall carry this out in the next section.

\section{Locally free sheaves on genus-zero curves}
\mylabel{Locally free}

We now make a digression and study locally free sheaves on regular genus-zero curves.
Recall that by Grothendieck's Splitting Theorem, every vector bundle over the Riemann sphere
is a sum of line bundles (\cite{Grothendieck 1957}, Theorem 2.1). 
This carries over to  the projective line over any ground field
(as in \cite{Okonek; Schneider; Spindler 1980}, Theorem 2.1.1 or 
\cite{Hazewinkel; Martin 1982}, Theorem 4.1).
One may rephrase the result by saying that the \emph{indecomposable} locally free sheaves
are exactly the invertible sheaves.
The situation for elliptic curves over algebraically closed fields is already much more complicated,
and was   solved by Atiyah \cite{Atiyah 1957}. The indecomposable sheaves over twisted lines
where determined by Biswas and Nagaray \cite{Biswas; Nagaraj 2009}, Theorem 4.1 
and  Novakovi\'c \cite{Novakovic 2012}, Corollary 6.1, using Galois descent. 

Fix a ground field $F$ of characteristic $p\geq 0$.
Let $D$ be a regular genus-zero curve  whose Picard group is generated by the dualizing sheaf,
that is, a twisted line or a twisted ribbon.
The dualizing  sheaf $\omega_D$ is invertible of degree $d=-2$, hence 
$\chi(\omega_D^{\otimes t}) = -2t +1$. Moreover,  
\begin{equation}
\label{omega powers}
h^1(\omega_D^{\otimes t}) =  h^0(\omega_D^{\otimes 1-t}) =0\quadand 
h^0(\omega_D^{\otimes t}) =  -2t +1\geq 1
\end{equation} 
provided $t\leq 0$ holds. 
The vector space $\Ext^1(\O_D,\omega_D)=H^1(C,\omega_D)$ is 1-dimensional, and the resulting non-split extension 
\begin{equation}
\label{indecompsalbe rank two}
0\lra \omega_D\lra \shF_D\lra \O_D\lra 0
\end{equation} 
defines a locally free sheaf $\shF_D$ of rank two with $\det(\shF_D)=\omega_D$. 
Up to isomorphism, it does not depend on
the choice of the extension, and is canonically attached to the curve $D$. We remark in passing that
this  construction works on every Gorenstein curve without constant field extension.
As we shall see, it plays a particularly important role for genus-zero curves.

Note that in the long exact sequence for \eqref{indecompsalbe rank two}, the connecting map $H^0(D,\O_D)\ra H^1(D,\omega_D)$
is non-zero, because the extension is non-split, hence it is bijective, and thus  $h^0(\shF_D)=h^1(\shF_D)=0$.

\begin{proposition}
\mylabel{rank two}
The sheaf $\shF_D$ is indecomposable, with $\chi(\shF_D\otimes \omega_D^{\otimes t}) = -4t$. Moreover,
$$
h^1(\shF_D\otimes \omega_D^{\otimes t})    =0\quadand 
h^0(\shF_D\otimes \omega_D^{\otimes t}) =  -4t\geq 4
$$
provided   $t\leq -1$ holds.
\end{proposition}

\proof
Suppose the sheaf is decomposable, and write $\shF_D=\omega_D^{\otimes m-2}\oplus\omega_D^{\otimes -m}$ for some  integer $m$.
From $h^0(\shF_D)=0$ we see $m-2,-m>0$, which gives the contradiction $0>m>2$. Thus $\shF_D$ is indecomposable.
Recall that for locally free sheaves $\shE_i$ of rank $r_i$ and degree $d_i$ we have 
\begin{equation}
\label{degree and rr}
\deg(\shE_1\otimes\shE_2)=  r_2d_1+r_1d_2\quadand \chi(\shE_i)=d_i+r_i,
\end{equation} 
the latter by Riemann--Roch, and the formula 
for $\chi(\shF_D\otimes \omega_D^{\otimes t})$ follows. 
Suppose now $t\leq -1$. Then the group $H^1(D,\omega_D^{\otimes t+1})$ vanishes. So
tensoring \eqref{indecompsalbe rank two} with $\omega_D^{\otimes t} $ yields an exact sequence 
$$
0\lra H^0(D,\omega_D^{\otimes t+1})\lra H^0(D,\shF_D\otimes\omega_D^{\otimes t})\lra H^0(D, \omega_D^{\otimes t})\lra 0,
$$
and the formula for $h^0(\shF_D\otimes\omega_D^{\otimes t})$ results from \eqref{omega powers}.
Furthermore, we get an identification $H^1(D,\shF_D\otimes\omega_D^{\otimes t})=H^1(D, \omega_D^{\otimes t})$, 
which indeed vanishes.
\qed

\medskip
Note that  the wedge product $\shF_D\otimes\shF_D\ra\det(\shF_D)$ defines   identifications
$$
\shF_D=\shF_D^\vee\otimes\omega_D\quadand \shF_D^\vee=\shF_D\otimes\omega_D^{\otimes-1},
$$
which is a special case of the following general result:

\begin{theorem}
\mylabel{indecomposable on genus-zero}
Up to isomorphism, the indecomposable locally free sheaves on $D$ are the 
$\omega_D^{\otimes a}$ and $\shF_D\otimes\omega_D^{\otimes b}$,
with exponents $a,b\in\ZZ$.
\end{theorem}

\proof
We have to show that
each locally free sheaf $\shE$ of finite rank  decomposes into a sum where each summand has
the form $\omega_D^{\otimes a}$ or $\shF\otimes\omega_D^{\otimes b}$, with various exponents $a$ and $b$.
We proceed  by induction on the rank $r\geq 0$. The cases $r=0$ is trivial.
Suppose now $r\geq 1$, and that the assertion holds for $r-1$.
Replacing $\shE$ by  a suitable $\shE\otimes\omega_D^{\otimes t}$, we may assume that $h^0(\shE)\neq 0$ but $h^0(\shE\otimes \omega_D)=0$.
If follows that there is a non-zero $s:\O_D\ra\shE$, and the  resulting cokernel $\shE'$ must be  locally free. We thus can view $\shE$
as an extension 
\begin{equation}
\label{first short exact}
0\lra \O_D\lra \shE\lra \shE'\lra 0.
\end{equation} 
By  the induction hypothesis, $\shE'$ is isomorphic to some 
$$
\shF(a_1,\ldots,a_m\mid b_1,\ldots,b_n) = 
\left(\bigoplus_{i=1}^m \omega_D^{\otimes a_i}\right)\times \left(\bigoplus_{j=1}^n \shF_D\otimes\omega_D^{\otimes b_j}\right).
$$
Tensoring the short exact sequence \eqref{first short exact}    with $\omega_D$, we obtain an exact sequence
$$
H^0(D,\shE\otimes \omega_D)\lra H^0(D,\shF(a_1+1,\ldots \mid b_1+1,\ldots )) \lra H^1(D,\omega_D).
$$
The term on the left vanishes,   the term on the right is one-dimensional, so the term in the middle is at most one-dimensional.
Using Proposition \ref{rank two}, we infer that $b_1+1, \ldots, b_n+1\geq 0$.
Arguing in a similar way with \eqref{omega powers}, we see that $a_i+1\geq 1$ for all $1\leq i\leq m$, with one possible
exception $i=s$, which than must have  $a_s+1=0$.  Summing up, we may assume that
$$
b_1,\ldots,b_n, a_1\geq -1\quadand a_2,\ldots,a_m\geq 0.
$$
For all exponents $t\geq 0$ and $s\geq -1$ the extension groups
\begin{equation}
\label{ext vanishing}
\Ext^1(\omega_D^{\otimes t},\O_D)=H^1(D,\omega_D^{\otimes -t})\quadand 
\Ext^1(\shF_D\otimes\omega^{\otimes s},\O_D)=H^1(\shF_D\otimes\omega_D^{\otimes -1-s})
\end{equation} 
vanishes. Consequently, \eqref{first short exact} splits provided  $a_1\geq 0$ holds,
and  $\shE$ decomposes as desired.

It remains to treat the case $a_1=-1$. The pullback $\shE_0=\shE\times_{\shE'}\omega_D^{\otimes a_1}$ of  \eqref{first short exact} with respect to 
the inclusion map $\omega_D^{\otimes a_1}\subset\shE'$ sits in a short exact sequence
\begin{equation}
\label{second short exact}
0\lra \O_D\lra \shE_0\lra \omega_D^{\otimes-1}\lra 0.
\end{equation} 
Each  such extension is  either isomorphic to $\shF_D^\vee$ or $\O_D\oplus \omega_D^{\otimes-1}$, because the vector space 
$\Ext^1(\omega_D^{\otimes-1},\O_D)=H^1(D,\omega_D)$ is one-dimensional.
Using the Snake Lemma, we see that our sheaves form a short exact sequence
\begin{equation}
\label{third short exact}
0\lra\shE_0\lra \shE\lra \shE''\lra 0,
\end{equation} 
now with $\shE''=\shF(a_2,\ldots,a_m\mid b_1,\ldots,b_n)$. The short exact sequence \eqref{second short exact} 
yields an exact sequence
$$
\Ext^1(\shE'',\O_D)\lra \Ext^1(\shE'',\shE_0)\lra \Ext^1(\shE'',\omega_D^{\otimes-1}).
$$
The outer terms vanish, as one sees by using \eqref{ext vanishing} again. In turn, the extension \eqref{first short exact} splits, and 
$\shE$ decomposes as desired.
\qed

\begin{proposition}
The indecomposable sheaves that are globally generated are precisely the $\omega_D^{\otimes a}$ with $a\leq 0$,
and $\shF_D\otimes\omega_D^{\otimes b}$ with $b\leq -1$. In any case, we have  
$$
\deg(\omega_D^{\otimes a}) =-2a\quadand \deg(\shF_D\otimes\omega_D^{\otimes b}) =   -4b-2.
$$
\end{proposition}

\proof
The degrees immediately follow from \eqref{degree and rr}. Clearly, $\omega_D^{\otimes a}$ is globally generated if and only
if $a\leq 0$. If $b\leq -1$, the outer terms in the  short exact sequence 
$0\ra \omega_D^{\otimes b+1}\ra \shF_D\otimes\omega_D^{\otimes b} \ra\omega_D^{\otimes b}\ra 0$
are globally generated, and $H^1(D,\omega_D^{\otimes b+1})$ vanishes, hence the term in the middle is globally generated. 
Conversely, if  $\shF_D\otimes\omega_D^{\otimes b}$ is globally generated, so is the quotient $\omega_D^{\otimes b}$,
and hence $b\leq 0$. The case $b=0$ is impossible, since we already observed that $h^0(\shF_D)=0$.
\qed

\medskip
For each locally free sheaf $\shE$ of rank two on $D$ yields the regular surface $X=\PP(\shE)$.
Each section  is of the form $C=\PP(\omega_D^{\otimes a})$, for some short exact sequence
\begin{equation}
\label{extension from section}
0\lra \omega_D^{\otimes a'}\lra \shE\lra \omega_D^{\otimes a}\lra 0.
\end{equation} 
The self-intersection is $C^2=\deg(\omega_D^{\otimes a})-\deg(\omega_D^{\otimes a'})= -2(a-a')$,
as explained in \cite{Fanelli; Schroeer 2020b}, Lemma 6.1.

\begin{proposition}
\mylabel{projectivization}
Notation as above. If $\shE$ is indecomposable, the selfintersection satisfies $C^2\equiv 2$ modulo $4$.
If $\shE$ is decomposable and $C^2>0$, then $C\subset X$ is linearly equivalent to a   curve $C'$ that is reducible. 
\end{proposition}

\proof
Suppose $\shE=\shF_D\otimes\omega_D^{\otimes b}$ for some integer $b$. Taking degrees for  \eqref{extension from section} we get
$$
2\equiv -4b-2=\deg(\shE)=-2a-2a'\equiv -2(a-a')= C^2\quad \text{modulo 4},
$$
which gives the first assertion. If $\shE$ is decomposable we have $\shE=\omega_D^{\otimes r}\oplus\omega_D^{\otimes r'}$
with $r+r'=a+a'$. Without loss of generality $r'\leq r$. Then the section $R=\PP(\omega_D^{\otimes r})$
has $R^2=-2(r-r')\leq 0$.
Let $F\subset X$ be the fiber over a closed point $b\in B$ of degree two. Then $C$ is numerically equivalent to $R+mF$
for some integer $m$. From
$0\leq (C\cdot R) = R^2 + m(F\cdot R) = R^2 + 2m$
we get  $m\geq -R^2/2\geq 0$. The case $m=0$ is impossible, because $C^2>0$ and $R^2\leq 0$, thus $R+mF$ is a reducible curve.
It is linearly  equivalent to $C$, because $\Pic^\tau(X)=0$.
\qed

\medskip
We close this section with some further observations on sheaves over such   $D$.
If $D$ is a twisted line, we have $D\otimes K\simeq\PP^1_E$ for some separable quadratic field extension $F\subset K$,
with
$$
\omega_D\mid \PP^1_K=\O_{\PP^1_K}(-2)\quadand \shF_D\mid\PP^1_K=\O_{\PP^1_K}(-1)\oplus \O_{\PP^1_K}(-1).
$$
If $D$ is a twisted ribbon, we are in characteristic $p=2$, 
and   $(D\otimes E)_\red=\PP^1_{E}$ for some height-one extension $F\subset E$ of degree four.
Note that the projection $\PP^1_E\ra B$ again has degree two, but now
$$
\omega_D\mid\PP^1_E=\O_{\PP^1_E}(-1)\quadand \shF_D\mid\PP^1_E=\O_{\PP^1_E}(-1)\oplus \O_{\PP^1_E}.
$$
Also note that the sheaf of K\"ahler differentials $\Omega^1_{D/k}$ is locally free of rank two.
It turns out to be  decomposable:

\begin{proposition}
\mylabel{kaehler differentials split}
In the above situation, we have $\Omega_{D/k}^1\simeq\omega_D^{\otimes 2}\oplus\omega_D$.
\end{proposition}

\proof
The sheaf $\omega_D^{\otimes-1}$ is very ample and embeds $D$ into $\PP^2$ as a curve of degree two.
The canonical surjection $\Omega^1_{\PP^2/F}|B\ra\Omega^1_{D/F}$ of locally free sheaves of rank two 
must be bijective. Restricting the  Euler sequence $0\ra \Omega^1_{\PP^2/F}\ra\bigoplus_{i=0}^2\O_{\PP^2}(-1)\ra\O_{\PP^2}\ra 0$
to $D$ yields a short exact sequence
$$
0\lra \Omega^1_{D/F}\lra  \omega_D^{\oplus 3}\lra \O_D\lra 0.
$$
This shows  $\deg(\Omega^1_{D/F})=-6$ and $\Hom(\Omega^1_{D/F},\omega_D)\neq 0$. 
In light of Theorem \ref{indecomposable on genus-zero}, the only possibilities 
for $\Omega^1_{D/F}$ are $\shF_D\otimes\omega_D$ and $\omega_D^{\otimes 2}\oplus\omega_D$.
Note that both have $h^0=0$ and $h^1=4$, so the numerical invariants do not tell them apart.

Choose a field extension $F\subset E$ such that the base-change $D\otimes E$ becomes the split ribbon 
$$
\PP^1_E\oplus\O_{\PP^1_E}(-1)=\Spec E[u,\epsilon] \cup \Spec E[u^{-1},\epsilon u^{-1}],
$$
where $\epsilon^2=0$. Over these affine open sets, the differentials $e_1=du,e_2=d\epsilon$
and $e'_1=du^{-1}$, $e'_2=d\epsilon u^{-1}$ form a basis for $\Omega^1_{D/F}\otimes E$, and the resulting
cocycle   takes the form
$$
\begin{pmatrix}u^{-2}&\epsilon u^{-2}\\0& u^{-1}\end{pmatrix} 
\in \GL_2(E[u^{\pm 1}, \epsilon]),
$$
as explained in \cite{Schroeer 2007a}, Section 2. With the new basis $e'_1$, $e'_2+\epsilon e'_1$
we get
$$
e'_1=u^{-2} e_1\quadand  e'_2+\epsilon e'_1 = u^{-1}e_2,
$$
and infer that the sheaf $\Omega^1_{D/F}$ and $\omega_D^{\otimes 2}\oplus\omega_D$ become isomorphic after base-change   to $  E$.
In turn, they are already isomorphic over $F$.
\qed

\medskip
Let us finally point out that the classification of locally free sheaves 
on the split ribbon $D\otimes E=\PP^1_E\oplus\O_{\PP^1_E}(-1)$ is much more complicated.
This can already seen on  a base-change $D\otimes K$ that is a denormalization of $\PP^1_E$, obtained by
replacing an $E$-valued point $x\in\PP^1_E$ by a $K$-valued point. To give a locally free sheaf on this denormalization
amounts to give a sheaf 
$\shF=\O_{\PP^1_E}(a_1)\oplus\ldots\oplus \O_{\PP^1_E}(a_r)$ 
on $\PP^1_E$, 
together with an $K$-rational structure on  the $E$-vector space $V=\shF(a)$, that is, an $K$-subspace $V_0$ such that the
canonical map $V_0\otimes_KE\ra V$ is bijective. This apparently involves continuous families.
It would be interesting to work this out.

\section{Proofs for main results}
\mylabel{Proofs}

Recall that  $X$ denotes  a minimal regular proper surface over a ground field $F$ 
with invariants $h^0(\O_X)=1$ and $h^1(\O_X)=h^0(\omega_X^{\otimes 2})=0$.
We now collect our findings and give the proofs for two main results of this paper.

\medskip
\emph{Proof for Theorem \ref{not integral}.} 
Here we have an integral curve  $C\subset X$ with $h^1(\O_C)=0$ and $C^2=4$, and the task is 
to find a  linearly equivalent  curve $C'$ that is not integral.
For this we have introduced  in Section \ref{Bend and break} a finite morphism $f:X\ra \PP^3$.
Several easy situations where already treated in Proposition \ref{thm in easy cases}, and it remains to handle the case
that the image $V=f(X)$  is an integral quartic surface, which is   non-normal by Proposition \ref{non-normal}.
We write  $g:Y\ra V$ for the  normalization, and only have to deal with the case that no singular point $y\in Y$
maps to a rational point $v\in V$, again by Proposition   \ref{thm in easy cases}. 

In light of Proposition \ref{degree branch curve} and Proposition \ref{thm for exotic cubics},
it remains to treat  the case that the branch curve $B\subset V$ for the normalization
is a twisted cubic. We then established in Theorem \ref{structure of blowing-up} that $X$ arises as blowing-up of the quartic surface $V$ with center the
twisted cubic $B$. From this we inferred that 
$X=\PP(\shE)$ for some locally free sheaf $\shE$ of rank two over a regular genus-zero curve $D$, having $C$ as a section.
For $D=\PP^1$ the sheaf $\shE$ is decomposable, and we find a reducible $C'$ as in the proof for Proposition \ref{projectivization}.
The interesting case is when $\omega_D$ generates the Picard group.
Recall that $C^2=4$. By Proposition \ref{projectivization}, the sheaf $\shE$ must be decomposable,
and $C$ is linearly equivalent to a reducible curve.
\qed

\medskip
\emph{Proof for Theorem \ref{trichotomy}.} 
The task is to show that $X$ is isomorphic to a plane or a quadric surface in $\PP^3$, or there is a fibration $f:X\ra B$ with $\O_B=f_*(\O_X)$
where the base and the generic fiber are genus-zero curves, or the dualizing sheaf $\omega_X$ generates the Picard group $\Pic(X)$.
Suppose the dualizing sheaf does not generate.
According to Proposition \ref{good curve} there is an integral curve $C\subset X$ with  $h^1(\O_C)=0$ and $C^2\geq 0$, such that every linearly
equivalent curve $C'$ is integral. In light of Proposition \ref{small self-intersection} and Proposition \ref{selfintersection numbers}
it remains to treat the case where  $C^2\geq 3$ and $C^2\mid 4$, in other words $C^2=4$.
But then  Theorem \ref{not integral} tells us that some $C'$ is non-integral, contradiction.
\qed


\end{document}